\documentclass{article}
\usepackage{amsmath}
\usepackage{amssymb}
\usepackage{mathrsfs}
\usepackage{amsthm}
\usepackage{stmaryrd}
\usepackage{geometry}
\usepackage{mathtools}
\usepackage{biblatex}

\usepackage{hyperref}
\usepackage[nameinlink,noabbrev]{cleveref}

\hypersetup{
    colorlinks,
    citecolor=black,
    filecolor=black,
    linkcolor=black,
    urlcolor=black,
    breaklinks=true
}

\def\W{\mathcal W}
\def\R{\mathbb R}
\def\N{\mathbb N}
\def\P{\mathbb P}
\def\E{\mathbb E}
\def\eps{\varepsilon}
\def\p{\partial}
\def\d{\mathrm{d}}
\def\endproof{\hfill$\Box$}
\newcommand{\e}{\mathrm{e}}

\DeclareMathOperator{\supp}{supp}
\DeclareMathOperator{\var}{Var}
\DeclareMathOperator{\re}{Re}
\DeclareMathOperator{\im}{Im}

\DeclarePairedDelimiterX{\wick}[1]{:}{:}{#1}

\crefname{section}{Section}{Sections}

\newtheorem{theo}{Theorem}[section]
\crefname{theo}{Theorem}{Therorems}

\newtheorem{prop}[theo]{Proposition}
\crefname{prop}{Proposition}{Propositions}

\newtheorem{lemme}[theo]{Lemma}
\crefname{lemme}{Lemma}{Lemmas}

\newtheorem{cor}[theo]{Corollary}
\crefname{cor}{Corollary}{Corollaries}

\theoremstyle{definition}
\newtheorem*{demo}{Proof}
\newtheorem{rem}[theo]{Remark}

\bibliography{biblio.bib}

\title{Wellposedness of the cubic Gross-Pitaevskii equation with spatial white noise on $\mathbb{R}^2$}

\author{
    P. Mackowiak\\ 
    CMAP, CNRS, École polytechnique,\\ 
    Institut Polytechnique de Paris, 91120 Palaiseau, France\\ 
    \texttt{pierre.mackowiak@polytechnique.edu}
}

\begin{document}

\maketitle

\begin{abstract}
    In this paper, we prove the global wellposedness of the Gross-Pitaevskii equation with white noise potential, i.e. a cubic nonlinear Schrödinger equation with harmonic confining potential and spatial white noise multiplicative term. This problem is ill-defined and a Wick renormalization is needed in order to give a meaning to solutions. In order to do this, we introduce a change of variables which transforms the original equation into one with less irregular terms. We construct a solution as a limit of solutions of the same equation but with a regularized noise. This convergence is shown by interpolating between a diverging bound in a high regularity Hermite-Sobolev space and a Cauchy estimate in $\mathbb{L}^2(\mathbb{\R}^2)$.\\
\end{abstract}

\noindent AMS: 35Q55, 60H40, 35R60

\noindent Keywords: Dispersive PDEs, Nonlinear Schrödinger equation, White noise

\section{Introduction}

In this paper, we are interested in solving the renormalized cubic Gross-Pitaevskii equation with spatial white noise potential in two dimensions of space
\begin{equation}{~}
    \begin{cases}\mathrm{i}\p_t u+ H u + \xi u+\lambda \left|u\right|^{2}u=0\\ u(0)=u_0,\end{cases} \label{Eq:2D}
\end{equation} 
where $H=\Delta-|x|^2$ is the Hermite operator. We do not expect, in general, \cref{Eq:2D} to have solutions, even for smooth and localized initial datum, due to the low regularity of the spatial white noise. 
Instead, we are interested in solving the equation 
\begin{equation}{~}
    \begin{cases}\mathrm{i}\p_t u+ \wick{H+\xi}u+\lambda \left|u\right|^{2}u=0\\ u(0)=u_0,\end{cases} \label{Eq:2D-renorm}
\end{equation} 
where $\wick{H+\xi}$ is the renormalized Hermite operator with spatial white noise potential, $\lambda\in\R$ and $u_0$ is in a certain function space which will be defined later. The unknown $u$ is a random complex-valued function defined on $I\times \mathbb{R}^2$, where $I$ is an interval of $\R$ containing 0.\\ 

This renormalized equation is defined through the same kind of exponential transformation as the one used in \cite{Haire_Labbe} for the continuous parabolic Anderson model, and in \cite{debussche_weber_T2} for the nonlinear Schrödinger equation on the 2d torus. Let $Y=(-H)^{-1}\xi$. Formally, the change of variable $u=v \e^Y$ yields $u$ is solution to \cref{Eq:2D} if and only if $v$ is solution to the transformed equation 
\begin{equation}{~}
    \begin{cases}\mathrm{i}\p_t v+ Hv +2\nabla Y\cdot \nabla v+ xY\cdot xv + |\nabla Y|^2 v +\lambda \left|v\right|^{2}v\e^{2 Y}=0\\ v(0)=u_0\e^{-Y}=v_0.\end{cases} \label{Eq:2D-transformed}.
\end{equation}
Since $\xi$ is, a priori, only in $\W^{-1-s,q}$ a.s., with $q>2$, $s>0$ and $sq>2$ (see \cref{Property_noise_regularity}), we only have $Y\in\W^{1-s,q}$ a.s., and thus $\nabla Y\in\W^{-s,q}$ and the term $|\nabla Y|^2$ is ill-defined. Here, spaces $\W^{s,q}$ denote Sobolev spaces associated with $-H$. We give more details on these spaces in \cref{subsection_sobolev_spaces}. In order to give a meaning to this equation, we have to renormalize the term $|\nabla Y|^2$ using non-homogeneous Wick products associated to $-H$, which were introduced in \cite{debouard_2d_GP}. We detail this renormalization in \cref{renorm_subsec}. We say that $u$ is solution to \cref{Eq:2D-renorm} if $u=v \e^Y$ and $v$ is solution to the renormalized transformed equation 
\begin{equation}{~}
    \begin{cases}\mathrm{i}\p_t v+ Hv +2\nabla Y\cdot \nabla v+ xY\cdot xv + \wick{|\nabla Y|^2} v +\lambda \left|v\right|^{2}v\e^{2 Y}=0\\ v(0)=v_0,\end{cases} \label{Eq:2D-transformed-renorm}
\end{equation}
where $\wick{|\nabla Y|^2}$ is defined as the limit in a space of negative regularity of
$$\wick{|\nabla Y_N|^2} = |\nabla Y_N|^2-\E\left[|\nabla Y_N|^2\right]$$
and $Y_N$ is a suitable regularization of $Y$ defined in \cref{subsection_smooth_truncation}. We emphasize that \cref{Eq:2D-transformed-renorm} has, formally, two invariant quantities, a transformed mass for $v\in\mathbb{L}^2(\mathbb{R}^2)$
\begin{equation}
    \Tilde{M}(v) = \int_{\mathbb{R}^2} |v|^2 \e^{2 Y} \d x = |u|_{\mathbb{L}^2_x}^2 \label{Eq:def-trans-mass}
\end{equation}
and a transformed energy for $v\in\W^{1,2}$
\begin{equation}
    \Tilde{\mathcal{E}}(v) = \frac{1}{2} \int_{\mathbb{R}^2} \left(\left|\nabla v\right|^2 + |xv|^2 - |xv|^2 Y \right)\e^{2 Y} \d x - \frac{1}{2}\left\langle \wick{|\nabla Y|^2}, |v|^2 \e^{2Y} \right\rangle - \frac{\lambda}{4} \int_{\mathbb{R}^2} |v|^4 \e^{4 Y}\d x, \label{Eq:def-trans-energy}
\end{equation}
where the bracket $\left\langle \wick{|\nabla Y|^2}, |v|^2 \e^{2Y} \right\rangle$ is defined as a duality bracket between suitable spaces.\\ 

This equation is in fact a nonlinear Schrödinger equation (NLS) on $\mathbb{R}^2$ with harmonic confining potential and spatial white noise potential. The equation without the harmonic potential was first studied by A. Debussche and H. Weber in \cite{debussche_weber_T2} on the 2d torus, where the authors proved the wellposedness of a renormalized problem for cubic nonlinearities and well-chosen random initial data. In order to prove this result, they solved the renormalized equation for $v_\eps = u_\eps \e^{-\Delta^{-1}\xi_\eps}$, where $\xi_\eps$ is a mollification of the noise, and then passed to the limit in probability as $\eps$ goes to 0. In order to do this, they showed a diverging $H^2$ bound with logarithmic loss in $\eps$ and a converging $\mathbb{L}^2(\mathbb{T}^2)$ bound. This wellposedness result was extended to fourth order nonlinearities by N. Tzvetkov and N. Visciglia in \cite{tzvetkov2020dimensional} using modified energies. They further extended the result to any positive exponent in \cite{Tzvetkov_2023} using Strichartz estimates with losses for $-\Delta+\xi_\eps$ in addition to modified energies.\\ 

The case of the full space, still without confining potential, was then studied in \cite{debussche2017solution} for sub-cubic nonlinearities, where they used weighted Besov spaces in order to deal with the growth of the white noise at infinity. It was then extended in \cite{debussche2023global} to any polynomial nonlinearities using Strichartz estimates with losses for $-\Delta+\xi_\eps$ on $\mathbb{R}^2$ in weighted Besov spaces. A similar result has been recently obtained for logarithmic nonlinearity in \cite{chauleur2023logarithmic}, also using weighted spaces.\\ 

Another approach is to work directly on the Anderson operator $-\Delta+\xi$, which has been widely studied recently, both on compact 2d manifolds (see \cite{allez2015continuous,bailleul2022analysis,gubinelli2018semilinear,Labb__2019,mouzard2021weyl,mouzard2023simple}) and on the full space (see \cite{ugurcan2022anderson}). Properties such as discrete spectrum, spectral gap or estimation of the eigenvalues have been obtained. Moreover some Strichartz estimates for this operator allow to solve in a more standard way NLS equations (see for example \cite{mouzard2022strichartz,zachhuber2020strichartz}). Up to our knowledge, there is no comparable result for the Hermite-Anderson operator $H+\xi$ yet, hence we chose not to follow this approach.\\ 

On the other hand, the deterministic equation with confining potential is an extensively studied dispersive equation known as the Gross-Pitaevskii equation (GP). It appears in the study of Bose-Einstein condensates as a mean-field limit to the coupled linear many-body Schrödinger equation with confining potential (see \cite{BEC_Pitaevskii}). The deterministic Gross-Pitaevskii equation is known to have different properties than the nonlinear Schrödinger equation due to the confining potential. For example, the harmonic potential allows standing waves to exists in the defocusing case (see \cite{ROSE1988207}) whereas such solutions cannot exist without the potential. Moreover, in dimension 1, it was shown in \cite{robert2022focusing} that the Gibbs measures with mass cut-off for the deterministic Gross-Pitaevskii equation are never normalizable at the critical exponent. This contrasts with the phase transition existing for the nonlinear Schrödinger equation on bounded domain. Up to our knowledge, there is no comparable result in 2d. A Gibbs measure which, formally, is invariant for the GP equation has been constructed in \cite{debouard_2d_GP}, but it was not proven to be invariant under the GP flow.\\ 

The two equations are linked by the so-called lens-transform that allows to pass from one to the other. This result can be used to prove various results for the NLS equation from the GP equation (see for example \cite{burq2020sure,Niederer1973TheMK,tao2009pseudoconformal,thomann2009random}). However, in presence of a white noise potential $\xi$, the lens-transform changes the white noise into a time dependent stochastic process with white noise distribution at each fixed time. This time dependency makes the lens-transform unsuitable for our purpose.\\ 

The Gross-Pitaevskii equation with random potential can be seen as a model for Bose-Einstein condensates in presence of inhomogeneities (see for example \cite{PhysRevLett.104.193901,PhysRevLett.109.243902} and references therein). Hence the white noise potential can be seen as a toy model of random potential with null correlation length. However, we do not know exactly what physical meaning should the Wick renormalization have in this case.\\ 

Mathematically, this case is close to the periodic NLS with white noise potential, as $H$ is a selfadjoint operator with compact resolvent on $\mathbb{R}^2$, as the Laplacian is on the torus. In this paper, inspired by the periodic case, we first prove the global wellposedness of the linear equation 
\begin{equation}
    \begin{cases}\mathrm{i}\p_t v+ Hv +2\nabla Y\cdot \nabla v+ xY\cdot x v + \wick{|\nabla Y|^2} v =0\\ v(0)=v_0\end{cases} \label{Eq:2D-transformed-renorm-linear}
\end{equation}
in some Bochner spaces, and then the global wellposedness of the cubic case (i.e. \cref{Eq:2D-transformed-renorm}) for any well-chosen random initial data, in both linear and cubic cases, provided the $\mathbb{L}^2_x$ norm of the initial data is small in the focusing case (i.e. $\lambda>0$). More precisely, we show the following results.

\begin{theo}{Linear case}\label{Theo-lin}

    Let $\sigma\in(1,2)$, $p_0\in(2,+\infty)$, $p\in(2,p_0)$, $a<0<b$ and $v_0\in\mathbb{L}^{p_0}(\Omega,\W^{2,2})$. The following hold:
    \begin{itemize}
        \item There exists a unique solution $v$ to \cref{Eq:2D-transformed-renorm-linear} in $\mathbb{L}^p\left(\Omega,\mathcal{C}([a,b],\W^{\sigma,2})\right)$
        \item The sequence $\left(v^{N}\right)_{N\in\N}$ of solutions of the regularized linear equation (\ref{Eq:2D-transformed-renorm-regularized-linear}) converges to $v$ in $\mathbb{L}^p\left(\Omega,\mathcal{C}([a,b],\W^{\sigma,2})\right)$
        \item The transformed mass and energy, defined by \cref{Eq:def-trans-mass,,Eq:def-trans-energy} with $\lambda=0$ respectively, are conserved.
    \end{itemize}
\end{theo}

\begin{theo}{Cubic defocusing case}\label{Theorem-cubic-defocusing}

    Let $\lambda\leqslant0$. Almost surely, for every $a<0<b$, every $\sigma\in(1,2)$ and every $v_0\in\W^{2,2}$, it holds
    \begin{itemize}
        \item there exists a unique solution $v$ to \cref{Eq:2D-transformed-renorm} in $\mathcal{C}([a,b],\W^{\sigma,2})$
        \item the sequence $\left(v^{N}\right)_{N\in\N}$ of solution of the regularized cubic equation (\ref{Eq:2D-transformed-renorm-regularized}) converges to $v$ in $\mathcal{C}([a,b],\W^{\sigma,2})$
        \item transformed mass and energy, defined by \cref{Eq:def-trans-mass,,Eq:def-trans-energy} respectively, are conserved.
    \end{itemize}
\end{theo}

\begin{theo}{Cubic focusing case}\label{Theorem-cubic-focusing}

    Let $\lambda>0$ and $L>0$, and let $\overline{\Omega}_L=\left\{\lambda L^2 \left|\e^{-2 Y}\right|_{\mathbb{L}^\infty_{x}}^2\left|\e^{2 Y}\right|_{\mathbb{L}^\infty_{x}}\left|\e^{4 Y}\right|_{\mathbb{L}^\infty_{x}}<4\right\}$. Then, on $\overline{\Omega}_L$, endowed with the conditional probability and $\sigma$-algebra, there exists a random variable $N_0\in\N$ a.s. such that almost surely, for every $a<0<b$, every $\sigma\in(1,2)$ and every $v_0\in\W^{2,2}$ verifying $|v_0|_{\mathbb{L}^2_{x}}\leqslant L$, it holds
    \begin{itemize}
        \item there exists a unique solution $v$ to \cref{Eq:2D-transformed-renorm} in $\mathcal{C}([a,b],\W^{\sigma,2})$
        \item the sequence $\left(v^{N}\right)_{N\geqslant N_0}$ of solution of the regularized cubic equation (\ref{Eq:2D-transformed-renorm-regularized}) converges to $v$ in $\mathcal{C}([a,b],\W^{\sigma,2})$
        \item transformed mass and energy, defined by \cref{Eq:def-trans-mass,,Eq:def-trans-energy} respectively, are conserved.\\ 
    \end{itemize}
\end{theo}

In both linear and cubic cases, we will follow the main ideas of \cite{debussche_weber_T2}, with some technical modifications. First, we regularize the noise using the smooth spectral truncation introduced in \cite{debouard_2d_GP} and detailed in \cref{subsection_smooth_truncation} instead of using a mollification. This is a more natural regularization in this context because it has a better behaviour in Sobolev-Hermite spaces than the usual truncation and regularization procedure. Then, we show a diverging bound in sufficiently regular space. Contrary to the periodic case, we only have a polynomial (and not a logarithmic) loss, but with an exponent as small as desired. We can obtain this diverging bound in $\W^{2,2}$ in the linear case but, in the nonlinear case, due to technical issues in the proof of \cref{Prop:diverging-bound-W^(sigma;2)-defocusing}, we only show a diverging bound in $\W^{\sigma,2}$ for any $\sigma\in(1,2)$. Finally, we obtain a converging $\mathbb{L}^2_x$ bound for the difference of approximations, and we use it along with the diverging bound and interpolation to prove the $\mathbb{L}^p_\omega$ or pathwise (depending on the case) convergence along a subsequence. The convergence along the whole sequence and conservation laws then follow easily.\\ 

Let us emphasize the differences with the periodic case treated in \cite{debussche_weber_T2}. First, notice that the quantity we subtract in order to renormalize $|\nabla Y_N|^2$ is a function of the space variable going to infinity as $N$ goes to infinity. In the periodic case, this quantity is constant in space, so the renormalization can be understood as an "infinite phase shift". This interpretation is no longer possible in presence of the harmonic potential. The renormalization should be thought as an "infinite potential" which prevents the solution from having infinite energy due to the high modes. Moreover, the eigenvalues of the Hermite operator grow slower than those of the periodic Laplacian. Thus, it is not sure that the diverging quantity goes to infinity as $\log(N)$ as it is the case in the periodic setting. In fact, one can show that it may diverge slower than $N^\eps$ for any $\eps>0$, but obtaining an optimal uniform upper bound on the divergence is still an open question. This leads to polynomial losses in \cref{Prop:W22-bound,norm-L2-diff-NM,Prop:diverging-bound-W^(sigma;2)-defocusing,Proposition-Bound-L2-Diff-cubic-defocusing} instead of the logarithmic loss in the periodic case. For this reason, it is necessary to take care not to have bounds of the form $N^{C \kappa}$ with a random variable $C$ depending on $\kappa$.\\ 

The paper is structured as follows, the \cref{subsection_hermite_op} recalls most results we need about the Hermite operator, the functional spaces we need in our analysis and their properties are introduced in \cref{subsection_sobolev_spaces}, the smooth spectral truncation is detailed in \cref{subsection_smooth_truncation}. The regularity properties of the white noise and other random variables involved in our analysis are presented in \cref{subsection_noise_regularity}, the renormalization procedure is explained in \cref{renorm_subsec} and the global wellposedness of the regularized equation is shown in \cref{subsection_regu_eq}. Finally, \cref{Section_linear_case} is devoted to the proof of \cref{Theo-lin} and \cref{Section_nonlinear_case} is devoted to the proof of \cref{Theorem-cubic-defocusing}. Most of the proofs of \cref{section_preliminaries,Section_noise} can be found in \cref{appendix_proof_preli,appendix_proof_proba} respectively.

\section{Preliminaries}\label{section_preliminaries}

\subsection{Notations and conventions}

\begin{itemize}
    \item $\Delta = \p_{x_1}^2+\p_{x_2}^2$ is the Laplacian
    \item For $x\in\mathbb{R}^2$, $|x|^2=x_1^2+x_2^2$ and $\langle x\rangle^2 = 1+|x|^2$
    \item $\N=\{0,1,2,\dotsc\}$ and $\N^*=\{1,2,\dotsc\}$
    \item For $\alpha\in\N^2$, $|\alpha|=\alpha_1+\alpha_2$ and $\p^\alpha = \p_{x_1}^{\alpha_1}\p_{x_2}^{\alpha_2}$
    \item For $s\in\R$ and $p\in(1,+\infty)$, $W^{s,p}=\left\{u \in\mathcal{S}'(\mathbb{R}^2), \mathcal{F}^{-1}\left(\langle \eta \rangle^s \mathcal{F}u\right)\in\mathbb{L}^p(\mathbb{R}^2)\right\}$ denotes the Bessel potential space endowed with the norm $|u|_{W^{s,p}_x} =  |\mathcal{F}^{-1}\left(\langle \eta \rangle^s \mathcal{F}u\right)|_{\mathbb{L}^p_x}$, where $\mathcal{F}$ denotes the usual Fourier transform. These spaces extend the usual Sobolev spaces by complex interpolation.
    \item We denote by $\mathcal{S}(\R^d)$ the space of smooth rapidly decaying functions on $\R^d$ and by $\mathcal{S}'(\R^d)$ its continuous dual, the space of tempered distributions on $\R^d$. We endow the space of tempered distributions with the real bracket
    $$\forall T\in\mathcal{S}'(\mathbb{R}^2)\cap\mathbb{L}^1_{loc}(\mathbb{R}^2),\forall \phi\in\mathcal{S}(\mathbb{R}^2), \langle T,\phi\rangle_{\mathcal{S}',\mathcal{S}} = \re\left(\int_{\mathbb{R}^2} T(x) \overline{\phi(x)}\d x\right)$$
    so that for $T\in\mathbb{L}^2(\mathbb{R}^2)$, $\langle T,\phi\rangle_{\mathcal{S}',\mathcal{S}}=\langle T,\phi\rangle_{\mathbb{L}^2_x}$.
    \item We say that a sequence of random variables $(X_N)_{N\in\N}$ defined on a probability space $(\Omega,\mathcal{A},\P)$ verifies the \hypertarget{property_star}{Property (*)} if
    \begin{itemize}
        \item For almost all $\omega\in\Omega$, $\sup_{N\in\N}X_N(\omega)<+\infty$ and for all $N\in\N$, $X_N(\omega)>0$.
        \item For all $p\in[1,+\infty)$, $\sup_{N\in\N}\E\left[X_N^p\right]<+\infty$.
    \end{itemize}
\end{itemize}

\subsection{Hermite operator}\label{subsection_hermite_op}

We denote by $H=\Delta-|x|^2$ the Hermite operator. It is well-known that this operator has eigenfunctions $(h_k)_{k\in\N^2}$ which are given by $h_k = \psi_{k_1}\otimes\psi_{k_2}$ ($k=(k_1,k_2)$), where we denote by $\psi_{k_1}$ the $k_1$-th Hermite function. For any $k\in\N^2$, $h_k$ verifies the relation $-H h_k = \lambda_k^2 h_k$ where $\lambda_k^2=2|k|+2$. Moreover, $(h_k)$ is a complete orthonormal system of $\mathbb{L}^2(\mathbb{R}^2,\R)$.\\ 

The Schwartz space $\mathcal{S}(\mathbb{R}^2,\R)$ and its dual $\mathcal{S}'(\mathbb{R}^2,\R)$ can be characterized using brackets against Hermite functions. The following proposition is a direct corollary of Theorem 1 in \cite{HermiteExpansionsSimon}.

\begin{prop}
    If $T\in\mathcal{S}'(\mathbb{R}^2,\R)$, then $\sum_{|k|\leqslant N} \langle T,h_k\rangle_{\mathcal{S}',\mathcal{S}}\,h_k$ converges weakly to $T$ in $\mathcal{S}'(\mathbb{R}^2,\R)$.
\end{prop}

Hence, for any $\alpha\in\R$, one can define $(-H)^\alpha$ on $\mathcal{S}'(\mathbb{R}^2,\R)$ by
$$\forall T\in\mathcal{S}'(\mathbb{R}^2,\R),\; (-H)^{\alpha}T = \sum_{k\in\N^2} \lambda_k^{2\alpha}\langle T,h_k\rangle_{\mathcal{S}',\mathcal{S}}\, h_k$$
defined as the weak limit of $\sum_{|k|\leqslant N} \lambda_k^{2\alpha}\langle T,h_k\rangle_{\mathcal{S}',\mathcal{S}}\, h_k$. For complex valued tempered distributions, analogous results hold for expansions in terms of $\{h_k,\mathrm{i} h_k\}_{k\in\N^2}$. It follows that $(-H)^{-1}$ has a kernel given by 
\begin{equation}
    K(x,y)=\displaystyle\sum_{k\in\N^2}\frac{h_k(x)h_k(y)}{\lambda_k^2} \label{noyau_K}
\end{equation}
in the sense that for any $\phi\in\mathcal{S}(\mathbb{R}^2)$, we have
$$\forall x\in\mathbb{R}^2,\;(-H)^{-1}\phi(x)=\int_{\mathbb{R}^2} K(x,y) \phi(y)\d y.$$
This kernel has the following integrability properties.

\begin{prop}{(See proposition 5 in \cite{debouard_2d_GP})}\label{regu_kernel_-H-1} For any $n\in\N^*$, $r\geqslant 2$ and $\alpha<1-\frac{2}{r}$, $K^n\in\mathbb{L}^r_x\W^{\alpha,2}_y$.

\end{prop}

This proposition will be useful to prove the convergence of the Wick products and to prove \cref{Property_noise_regularity} about the white noise regularity.

\subsection{Hermite-Sobolev spaces}\label{subsection_sobolev_spaces}

Let $s\in\R$ and $p\in(1,+\infty)$ and define
$$\W^{s,p} =\left\{u\in\mathcal{S}'(\mathbb{R}^2), (-H)^{\frac{s}{2}}u\in\mathbb{L}^p(\mathbb{R}^2)\right\}$$
endowed with the norm $$|u|_{\W^{s,p}_x} = |(-H)^{\frac{s}{2}}u|_{\mathbb{L}^p_x}.$$
It is known that these spaces are Banach spaces and that for $q$ such that $\frac{1}{p}+\frac{1}{q}=1$, we have $\left(\W^{s,p}\right)'=\W^{-s,q}$ with equal norm. The subspace $span\left\{h_k, k\in\N^2\right\}$ is dense in these spaces (see for example \cite{BongioanniHSspaces} for the case $s\geqslant0$) and for $s\geqslant0$, an equivalent norm is given for $p\in(1,+\infty)$ by (see \cite{Dziubanski})
$$|\langle D\rangle^s u|_{\mathbb{L}^{p}_x}+|\langle x\rangle^s u|_{\mathbb{L}^p_x}.$$
This shows that $\W^{s,p}=\left\{u\in W^{s,p}(\mathbb{R}^2), \langle x\rangle^s u\in\mathbb{L}^p(\mathbb{R}^2)\right\}.$
Using this equivalent norm, we define $\W^{k,p}$ for $k\in\N$ and $p\in\{1,+\infty\}$ as
$$\W^{k,p} = \left\{u\in W^{k,p}(\mathbb{R}^2), \langle x\rangle^k u\in\mathbb{L}^p(\mathbb{R}^2)\right\}$$
endowed with the norm
$$|u|_{\W^{k,p}_x} = |u|_{W^{k,p}_x}+|\langle x\rangle^k u|_{\mathbb{L}^p_{x}}$$
and extend this definition to $s\geqslant 0$ by complex interpolation between integer regularities. Moreover inspired by the relation $\left(\W^{s,p}\right)'=\W^{-s,q}$, for $s<0$, we define $\W^{s,\infty} = \left(\W^{-s,1}\right)'$. As usual, we cannot define spaces of negative $\mathbb{L}^1$ regularity as dual spaces of positive $\mathbb{L}^\infty$ regularity.

\begin{rem}{~}
    \begin{itemize}
        \item We do not claim that the first definition of $\W^{s,p}$ agrees with the one we give for spaces over $\mathbb{L}^1$ and $\mathbb{L}^\infty$. We use the same notation for simplicity as we will always control $\W^{s,\infty}$ norms by Sobolev embeddings.
        \item In the special case $\W^{1,2}$, we also use the notation $\Sigma$.
    \end{itemize}
\end{rem}

The choices we made ensure us that we can interpolate between spaces of positive regularity and some other convenient properties that we need in our analysis.\\ 

A useful property of these spaces is that differentiation and multiplication by $x$ act directly on regularity. In order to see this, let us introduce the annihilation and creation operators:
$$\forall i\in\{1,2\},\; A_i = \p_i + x_i \text{ and } A_{-i} = A_i^* = -\p_i + x_i.$$

\begin{prop}\label{Proposition_creation_annhil_Wsp}

    Let $s\in\R$, $p\in(1,+\infty)$ and $i\in\{\pm 1,\pm 2\}$, then $A_i$ is a bounded linear operator from $\W^{s,p}$ to $\W^{s-1,p}$.
    
\end{prop}

\begin{cor}\label{Corollary_action_D_and_x_on_Wsp}\label{action_D&x/Wsp}

    Let $s\in\R$, $p\in(1,+\infty)$ and $i\in\{ 1,2\}$, then $\p_i$ and $x_i\cdot$ are bounded linear operators from $\W^{s,p}$ to $\W^{s-1,p}$.
    
\end{cor}

The Sobolev-Hermite spaces verify some continous embeddings as in  the classical Sobolev framework. We will also refer to these embeddings as Sobolev embeddings.

\begin{prop}{(Sobolev embeddings)}\label{Sobolev-embeddings}
    Let $1< p\leqslant q<+\infty$ and $s>\sigma$ such that $\frac{1}{p}-\frac{s}{2}\leqslant\frac{1}{q}-\frac{\sigma}{2}$. Then $\W^{s,p}$ is continuously embedded in $\W^{\sigma, q}$. Moreover if $1<p< q\leqslant+\infty$ and $s>\sigma$ such that $\frac{1}{p}-\frac{s}{2}<\frac{1}{q}-\frac{\sigma}{2}$, then $\W^{s,p}$ is compactly embedded in $\W^{\sigma,q}$.
    
\end{prop}

Moreover, our choices allow us to have a simple product rule on our spaces.

\begin{lemme}\label{product_rule}

    Let $s\geqslant0$ and $1\leqslant p,q,r\leqslant+\infty$ such that $\frac{1}{p}+\frac{1}{q}=\frac{1}{r}$. There exists a constant $C>0$ such that for all $u\in\W^{s,p}$ and $v\in\W^{s,q}$, we have $uv\in\W^{s,r}$ and $|uv|_{\W^{s,r}_x}\leqslant C |u|_{\W^{s,p}_x}|v|_{\W^{s,q}_x}$.
    
\end{lemme}

\begin{demo}

    For $s\in\N$, the claim follows from Leibniz's formula and Hölder's inequality. By bilinear interpolation (see for example theorem 4.4.1 of \cite{bergh2011interpolation}), we obtain the claim for all $s\geqslant 0$.\endproof
    
\end{demo}

We deduce the following product rule for negative-positive regularity products.

\begin{lemme}\label{product-rule-W(-k;r)-Wkp}
    Let $s\geqslant0$, $1<p<+\infty$ and $1\leqslant q,r\leqslant+\infty$ such that $1-\frac{1}{p}+\frac{1}{q}=1-\frac{1}{r}$. There exists a constant $C>0$ such that for all $u\in\W^{-s,r}$ and $v\in\W^{s,q}$, we have $uv\in\W^{-s,p}$ and $|uv|_{\W^{-s,p}_x}\leqslant C |u|_{\W^{-s,r}_x}|v|_{\W^{s,q}_x}$.
\end{lemme}

\begin{demo}

    Let $p'$ and $r'$ be such that $\frac{1}{p'}=1-\frac{1}{p}$ and $\frac{1}{r'}=1-\frac{1}{r}$. Let $\phi\in\mathcal{S}(\mathbb{R}^2)$, using \cref{product_rule} we have $\overline{v}\phi\in\W^{s,r'}$ whose dual is $\W^{-s,r}$. Then, we can define the tempered distribution $uv$ as $$\langle uv,\phi\rangle_{\mathcal{S}',\mathcal{S}} = \langle u,\overline{v}\phi\rangle_{\W^{-s,r}_x,\W^{s,r}_x} $$ and we have $|\langle uv,\phi\rangle_{\mathcal{S}',\mathcal{S}}|\lesssim |u|_{\W^{-s,r}_x}|\overline{v}\phi|_{\W^{s,r'}_x}\lesssim|u|_{\W^{-s,r}_x}|v|_{\W^{s,q}_x}|\phi|_{\W^{s,p'}_x}$. We conclude by density of Schwartz functions.\endproof
    
\end{demo}

Let us recall another product rule on these spaces whose proof is in \cite{debouard_2d_GP} and which will be usefull in the proof of \cref{conv_renorm}.

\begin{lemme}\label{product_rule_2var}

    Let $r,r_1,r_2> 1$, $p\geqslant 2$ and $0<\alpha<\sigma$ such that $\frac{1}{r_1}+\frac{1}{r_2}=\frac{1}{r}$ and $\frac{\sigma-\alpha}{2}=\frac{1}{p}$. There exists a constant $C>0$ such that for all $f,g\in\mathbb{L}^{r_1}_x\mathbb{L}^p_y\cap\mathbb{L}^{r_2}_x\W^{\sigma,2}_y$, it holds $fg\in\mathbb{L}^r_x\W^{\alpha,2}_y$ and $$|fg|_{\mathbb{L}^r_{x}\W^{\alpha,2}_{y}} \leqslant C\left(|f|_{\mathbb{L}^{r_1}_{x}\mathbb{L}^p_{y}}|g|_{\mathbb{L}^{r_2}_{x}\W^{\sigma,2}_{y}}+|f|_{\mathbb{L}^{r_2}_{x}\W^{\sigma,2}_{y}}|g|_{\mathbb{L}^{r_1}_{x}\mathbb{L}^p_{y}}\right).$$
     
\end{lemme}

\subsection{The smooth truncation}\label{subsection_smooth_truncation}

Choose $\chi\in\mathcal{S}(\R)$ such that $0\leqslant \chi\leqslant 1$, $\supp(\chi)\subset[-1,1]$ and $\chi_{|[0,1/2]}=1$. We define the smooth truncation as the family $(S_N)_{N\in\N}$ with $S_N= \chi\left(\frac{-H}{\lambda_N^2}\right)$ where $\lambda_N=\lambda_{(N,0)}=\sqrt{2N+2}$. For $u\in\mathbb{L}^{2}(\mathbb{R}^2)$, $S_N$ is well-defined and we have the formula 
$$ S_N u = \sum_{k\in\N^2} \chi_{k,N} \langle u,h_k\rangle_{\mathbb{L}^2}\, h_k,$$
where $\chi_{k,N}=\chi\left(\frac{\lambda_k^2}{\lambda_N^2}\right)$. It extends naturally to $\W^{s,2}$ with $s\geqslant 0$ and to the case $s<0$ by duality. Using Sobolev embeddings, one can extend this definition to any $\W^{s,p}$ with $s\in\R$ and $p\in(1,2)$ and finally to $\W^{s,q}$ with $s\in\R$ and $q\in(2,+\infty)$ by duality. One can show using Fourier-Hermite expansions that for $u\in\W^{s,p}$ ($s\in\R$ and $p\in(1,+\infty)$), we have
$$S_N u = \sum_{k\in\N} \chi_{k,N} \langle u, h_k\rangle_{\mathcal{S}',\mathcal{S}}\, h_k$$
which extends the previous $\mathbb{L}^2$ formula. We have the following important lemma which motivates the choice of a smooth truncation.

\begin{lemme}{(see Theorem 1 in \cite{Jensen1995LpmappingPO})}\label{Lemma_psi(-thetaH)}
Let $\psi\in\mathcal{S}(\R)$ and $p\in[1,+\infty]$. The family $(\psi(-\theta H))_{\theta\in[0,1]}$ is uniformly bounded in $\mathcal{L}\left(\mathbb{L}^p(\mathbb{R}^2)\right)$.
    
\end{lemme}

Using the fact that $S_N$ commutes with any $(-H)^{\alpha}$, we obtain the following result, which is not verified by the spectral projector $\Pi_N u = \sum_{|k|\leqslant N}\langle u,h_k\rangle_{\mathbb{L}^2}\, h_k$.

\begin{cor}\label{bound-S_N-Wsq}
    Let $s\in\R$ and $p\in(1,+\infty)$. Then, $(S_N)_{N\in\N}$ is uniformly bounded in $\mathcal{L}\left(\W^{s,p}\right)$.
    
\end{cor}
We will use the following lemmas, which are in fact corollaries of \cref{Lemma_psi(-thetaH)}, to estimate norms of low-mode and high-mode functions.
\begin{lemme}\label{lemme-low-freq-estim}
    Let $p\in(1,+\infty)$, $\alpha\in\R$ and $s>0$, there exists $C>0$ such that
    $$\forall N\in\N,\; \forall \phi\in\W^{\alpha,p},\; |S_N \phi|_{\W^{\alpha+s,p}_x}\leqslant C \lambda_N^{s}|\phi|_{\W^{\alpha,p}_x}.$$
\end{lemme}

\begin{lemme}\label{lemme-high-freq-estim-Ws2}
    Let $\alpha\in\R$ and $s>0$. We have
    $$\forall N\in\N,\; \forall \phi\in\W^{\alpha+s,2},\; |\phi-S_N \phi|_{\W^{\alpha,2}_x}\leqslant \lambda_{\left\lfloor\frac{N}{2}\right\rfloor}^{-s}|\phi|_{\W^{\alpha+s,2}_x}.$$
\end{lemme}

By density of Schwartz functions, Sobolev embeddings and \cref{bound-S_N-Wsq}, it follows that for any $s\in\R$, any $q\in[2,+\infty)$ and any $u\in\W^{s,q}$, $S_N u$ converges to $u$ in $\W^{s,q}$.

\section{The regularized noise}\label{Section_noise}
\subsection{Noise regularity}\label{subsection_noise_regularity}

Let $(\Omega,\mathcal{A},\P)$ be a probability space such that there exists a sequence $(\xi_k)_{k\in\N^2}$ of standard normal random variables. We define the spatial white noise on $\mathbb{R}^2$ as $\displaystyle\xi = \sum_{k\in\N^2} \xi_k h_k$. Such random Hermite series have been studied in \cite{imekraz2014random}. As $(h_k)_{k\in\N^2}$ is a complete orthonormal system of $\mathbb{L}^2(\mathbb{R}^2)$, one recovers the fundamental property of white noises
$$\forall\phi,\psi\in\mathcal{S}(\mathbb{R}^2),\; \E\left[\langle \xi,\phi\rangle_{\mathcal{S}',\mathcal{S}}\,\langle \xi,\psi\rangle_{\mathcal{S}',\mathcal{S}}\right] = \langle \phi,\psi\rangle_{\mathbb{L}^2_{x}}.$$

\begin{prop}\label{Property_noise_regularity}
    Let $q>2$ and $s>0$ such that $sq>2$. Then $\xi\in\W^{-1-s,q}$ a.s.
\end{prop}

\begin{demo}
     Let $\alpha>0$ and $q>2$, we have $(-H)^{-\alpha/2}\xi = \displaystyle\sum_{k\in\N^2} \xi_k \lambda_k^{-\alpha} h_k$. Using Fubini's theorem and gaussianity, we obtain 
     \begin{equation*}
         \E\left[\left|\xi\right|^{q}_{\W^{-\alpha,q}_x}\right] = \int_\R \E\left[\left|(-H)^{-\alpha/2}\xi(x)\right|^{q}\right] \d x \lesssim |\sigma_\alpha|_{\mathbb{L}^{q}_{x}}^{q}
     \end{equation*}
    where $$\sigma_\alpha(x)^2 = \displaystyle\sum_{k\in\N^2} \lambda_k^{-2\alpha}|h_k(x)|^2=\E\left[\left|(-H)^{-\alpha/2}\xi(x)\right|^{2}\right].$$
    Using the definition of the kernel $K$ in \cref{noyau_K}, we obtain $\sigma_\alpha(x)^2 =|K(x,\cdot)|_{\W^{2-\alpha,2}_y}^2$. Then by \cref{regu_kernel_-H-1}, $\sigma_\alpha$ has finite $\mathbb{L}^q_x$ norm when $2-\alpha<1-\frac{2}{q}$ (i.e. $\alpha>1+\frac{2}{q}$). Let $s=\alpha-1$, then $s>0$ and $\xi\in\W^{-1-s,q}$ a.s. as long as $qs>2$.\endproof
\end{demo}

\begin{rem}\begin{itemize}
    \item This proof is not totally rigorous, but in fact, the following lemmas will show the convergence of $S_N\xi$ to $\xi$ in $\W^{-1-s,q}$ ($q>2$, $sq>2$), both almost surely and in all $\mathbb{L}^p_\Omega$ spaces, with an explicit rate of convergence. Moreover, this regularity result is the same as the one given by applying Theorem 2.2 in \cite{imekraz2014random} to $(-H)^{-1-s}\xi$.
    \item We do not know if this regularity is optimal but we strongly believe that for $q<+\infty$, $\xi$ does not belong to $\W^{-1,q}$. This contrasts with the case of a radial white noise (i.e. with only radial Hermite functions in the expansion of $\xi$) which lives almost surely in $\W^{-1,q}$ for any $q\in(2,+\infty)$ according to Theorem 2.3 in \cite{imekraz2014random}.
\end{itemize}
    
\end{rem}

Now, we introduce a probabilistic lemma that we will use several times in what follows.

\begin{lemme}\label{Lemma_bound_from_moment}

    Let $(X_N)_{N\in\N}$ be a sequence of real random variables defined on the same probability space. Assume there exists $\alpha>0$ such that for every $p\in[1+\infty)$, $\displaystyle\sup_{N\in\N} \lambda_N^\alpha |X_N|_{\mathbb{L}^p_\omega} <+\infty.$
    Then, for every $\beta<\alpha$, almost surely $\displaystyle\sup_{N\in\N} \lambda_N^\beta |X_N| <+\infty.$
    
\end{lemme}

Let $N\in\N$ and set $Y=(-H)^{-1} \xi$ and $Y_N=S_N Y$. The proof of the following results are given in appendix.

\begin{lemme}{~}\label{Lem:convergence-Y_N-W(1-s;q)}
    Let $p\in[1,+\infty)$, $q>2$ and $1>s>s'>\frac{2}{q}$. There exists a constant $C>0$ such that 
    $$\sup_{N\in\N}\left[\lambda_N^{s-s'}|Y-Y_N|_{\mathbb{L}^p_\omega\W^{1-s,q}_x}\right]\leqslant C.$$

\end{lemme}

Using Sobolev embeddings and \cref{Lemma_bound_from_moment}, this implies the following result. 

\begin{cor}\label{Cor:convergence-Y_N-W(1-kappa;infty)}

    Let $0<\kappa'<\kappa<1$, there exists a sequence of positive random variables $(C_N)_{N\in\N}$ verifying \hyperlink{property_star}{Property (*)} such that almost surely, for all $N\in\N$,
    $$|Y-Y_N|_{\W^{1-\kappa,\infty}_x}\leqslant C_N \lambda_N^{-\kappa'},\; |\nabla Y-\nabla Y_N|_{\W^{-\kappa,\infty}_x}\leqslant C_N \lambda_N^{-\kappa'} \text{ and }|xY-xY_N|_{\W^{-\kappa,\infty}_x} \leqslant C_N \lambda_N^{-\kappa'}.$$

    Moreover, for any $p\in[1,+\infty)$, there exists $C>0$ such that for all $N\in\N$,
    $$|Y-Y_N|_{\mathbb{L}^p_\omega\W^{1-\kappa,\infty}_x}\leqslant C\lambda_N^{-\kappa'},\;|\nabla Y-\nabla Y_N|_{\mathbb{L}^p_\omega\W^{-\kappa,\infty}_x}\leqslant C\lambda_N^{-\kappa'} \text{ and }|xY-xY_N|_{\mathbb{L}^p_\omega\W^{-\kappa,\infty}_x}\leqslant C\lambda_N^{-\kappa'}.$$
\end{cor}

Using \cref{bound-S_N-Wsq,Property_noise_regularity,Sobolev-embeddings}, we obtain the following lemma.

\begin{cor}\label{Corollary_divergence_moment_norm_Linfty_GradYN_xYN}

    Let $\kappa\in(0,1)$, there exists a sequence of positive random variables $(C_N)_{N\in\N}$ verifying \hyperlink{property_star}{Property (*)} such that almost surely  
    $$\forall N\in\N,\; |\nabla Y_N|_{\mathbb{L}^{\infty}_{x}}\leqslant C_N\lambda_N^{\kappa} \text{ and }|x Y_N|_{\mathbb{L}^{\infty}_{x}}\leqslant C_N\lambda_N^{\kappa}.$$
    
    Moreover, for any $p\in[1,+\infty)$, there exists $C>0$ such that 
    $$\forall N\in\N,\; |\nabla Y_N|_{\mathbb{L}^p_\omega\mathbb{L}^{\infty}_{x}}\leqslant C\lambda_N^{\kappa} \text{ and }|x Y_N|_{\mathbb{L}^p_\omega\mathbb{L}^{\infty}_{x}}\leqslant C\lambda_N^{\kappa}.$$
    
\end{cor}

Using Hölder's inequality, Sobolev embeddings and \cref{Lemma_bound_from_moment}, we deduce the following result.

\begin{cor}\label{Cor:estimate-norm-Lq-X-Y_N}

    Let $\kappa\in(0,1)$, $q\in(1,+\infty)$ and a random variable $X\in\W^{1,q}$ a.s., there exists a sequence of positive random variables $(C_N)_{N\in\N}$ verifying \hyperlink{property_star}{Property (*)} such that almost surely 
    $$\forall N\in\N,\; |\nabla X\cdot\nabla Y_N|_{\mathbb{L}^{q}_x}\leqslant C_N |X|_{\W^{1,q}_x} \lambda_N^\kappa \text{ and }|x X\cdot x Y_N|_{\mathbb{L}^{q}_x}\leqslant C_N  |X|_{\W^{1,q}_x} \lambda_N^\kappa.$$
    
    Moreover, let $p_0>p\geqslant 1$ and assume that $X\in\mathbb{L}^{p_0}(\Omega,\W^{1,q})$, then there exists $C>0$ such that
    $$\forall N\in\N,\;|\nabla X\cdot\nabla Y_N|_{\mathbb{L}^p_\omega\mathbb{L}^{q}_x}\leqslant C \lambda_N^{\kappa} |X|_{\mathbb{L}^{p_0}_\omega\W^{1,q}}\text{ and }|x X\cdot x Y_N|_{\mathbb{L}^p_\omega\mathbb{L}^{q}_x}\leqslant C \lambda_N^{\kappa} |X|_{\mathbb{L}^{p_0}_\omega\W^{1,q}}.$$

\end{cor}

Using the algebra property of $\W^{1-s,q}$ for $q>4$ and $sq>2$, $\e^{a Y}-1\in\W^{1-s,q}$ for any $a\in\R$. Using the gaussianity of $Y_N$ and $Y$, we obtain the following moment estimates for exponentials of $Y_N$ and $Y$.

\begin{lemme}\label{Lem:moment-boud-expYN-Wsq}
    Let $p\geqslant 1$, $q>4$, $s\in\left(\frac{2}{q},1-\frac{2}{q}\right)$ and $a\in\R$, there exists $C>0$ such that 
    $$\left|\e^{a Y}-1\right|_{\mathbb{L}^p_\omega\W^{1-s,q}_x}\leqslant C \text{ and }\sup_{N\in\N}\left|\e^{a Y_N}-1\right|_{\mathbb{L}^p_\omega\W^{1-s,q}_x}\leqslant C.$$
\end{lemme}

Using Sobolev embeddings, we obtain the following corollary.

\begin{cor}\label{Cor:moment-boud-expYN-norm-Linf}
    Let $p\geqslant 1$, $\kappa\in(0,1)$ and $a\in\R$, there exists $C>0$ such that
    $$\left|\e^{a Y}-1\right|_{\mathbb{L}^p_\omega\W^{1-\kappa,\infty}_x}\leqslant C \text{ and } \sup_{N\in\N}\left|\e^{a Y_N}-1\right|_{\mathbb{L}^p_\omega\W^{1-\kappa,\infty}_x}\leqslant C.$$
\end{cor}

We deduce a control on exponential moments in $\mathbb{L}^\infty$ Hermite-Sobolev spaces.

\begin{lemme}\label{Lem:Conv-W(1-kappa;inf)-exp(aYN)}
    Let $0<\kappa'<\kappa<1$, $a\in\R$ and $p\geqslant 1$. There exists $C>0$ such that
    $$\sup_{N\in\N}\left[\lambda_N^{\kappa'}\left|\e^{aY_N}-\e^{aY}\right|_{\mathbb{L}^p_\omega\W^{1-\kappa,\infty}_x}\right]\leqslant C.$$

    Moreover, there exists a positive random variable $C$ almost surely finite such that, almost surely, it holds
    $$\forall N\in\N,\;  \left|\e^{aY_N}-\e^{aY}\right|_{\W^{1-\kappa,\infty}_x}\leqslant C\lambda_N^{-\kappa'} \text{ and } \left|\e^{a Y_N}-1\right|_{\W^{1-\kappa,\infty}_x}\leqslant C.$$
\end{lemme}

We introduce the shorten notation $|f|_{1+\W^{\sigma,p}_x}=1+|f -1|_{\W^{\sigma,p}_x}$ for any $f\in\mathcal{S}'(\mathbb{R}^2)$ verifying $f-1\in\W^{\sigma,p}$ ($\sigma\in\R$ and $p\in[1,+\infty]$). We emphasize that this is just a notation for the seek of compactness and that it does not denote the norm of any space. It allows to formulate in a compact fashion the following straightforward corollary of \cref{product_rule}.

\begin{lemme}\label{lem:product-with-exp(aY)}
    Let $s\in(0,1)$ and $p\in[1,+\infty]$. There exists a deterministic constant $C>0$ such that almost surely, for any $v\in\W^{s,p}$ and any $a\in\R$, $v \e^{a Y}\in\W^{s,p}$ and it holds
    $$\left|v \e^{a Y}\right|_{\W^{s,p}_x}\leqslant C \left|v \right|_{\W^{s,p}_x}\left| \e^{a Y}\right|_{1+\W^{s,\infty}_x}.$$
\end{lemme}

\subsection{Renormalization}\label{renorm_subsec}

Let $N\in\N^*$, define $\wick{|\nabla Y_N|^2}=|\nabla Y_N|^2-C_N^2$ where $C_N(x)^2=\E\left[|\nabla Y_N(x)|^2\right]=\displaystyle\sum_{k\in\N^2}\chi_{k,N}^2 \frac{|\nabla h_k(x)|^2}{\lambda_k^4}$. We now prove the following result on the convergence of $:|\nabla Y_N|^2:$.
\begin{prop}{~}\label{conv_renorm}
    Let $q>2$ and $s>\frac{2}{q}$. Then $(\wick{|\nabla Y_N|^2}:)_{N\in\N}$ is a Cauchy sequence in the Banach space $\mathbb{L}^q\left(\Omega,\W^{-s,q}\right)$, and thus converges to an element $\wick{|\nabla Y|^2}\in\mathbb{L}^q\left(\Omega,\W^{-s,q}\right)$. Moreover for all $0<\delta<\frac{2s}{3}-\frac{4}{3q}$ there exists a constant $C>0$ such that for all $M\in\N^*$
    $$\left|\wick{|\nabla Y_M|^2}-\wick{|\nabla Y|^2}\right|_{\mathbb{L}^q_\omega\W^{-s,q}_x}\leqslant C \lambda_M^{-\delta}.$$

\end{prop}

Let us recall some results about Wiener chaos and renormalization. For an introduction to Wick products see \cite{Prato2007WickPI} and for more details about Wick products adapted to $-H$, see \cite{debouard_2d_GP}. First recall the probabilist's Hermite polynomials:
$$\forall x\in\R,\; \forall n\in\N,\; H_n(x) = \frac{(-1)^n\e^{\frac{x^2}{2}}}{\sqrt{n!}}\frac{\d^n}{\d x^n}\left[e^{-\frac{x^2}{2}}\right].$$
They form a complete orthogonal system of $\mathbb{L}^2(\R,\mu)$ $\left(\text{ with }\d \mu(x) = \frac{\exp(-x^2/2)}{\sqrt{2\pi}}\d x\right)$.\\ 

Let $W:\mathbb{L}^2(\mathbb{R}^2)\to\mathbb{L}^2(\Omega)$ be the unique isometry such that $W_{h_k}=\xi_k$ for all $k\in\N^2$. We define the $n$-th chaos as 
$$ \mathcal{H}_n = \overline{\left\{H_n(W_f), |f|_{\mathbb{L}^2_{x}}=1\right\}}^{\mathbb{L}^2(\Omega)}.$$
It is well known that $\mathbb{L}^2(\Omega,\mathcal{G}) = \bigoplus_{n\in\N} \mathcal{H}_n$, where $\mathcal{G}=\sigma\left(\{W_f, f\in\mathbb{L}^2(\mathbb{R}^2)\}\right)=\sigma\left(\{\xi_k, k\in\N^2\}\right)$ is the $\sigma$-field generated by the white noise $W$. Thanks to the orthonormality of $(H_n)_{n\in\N}$, we have the following relation.

\begin{lemme}\label{esp-produi-W_fW_g}

    Let $k,n\in\N$ and $f,g\in\mathbb{L}^2(\mathbb{R}^2)$ such that $|f|_{\mathbb{L}^2_{x}}=1=|g|_{\mathbb{L}^2_{x}}$. Then, we have $\E\left[H_k(W_f)H_n(W_g)\right] = k! \langle f,g\rangle_{\mathbb{L}^2_{x}}^k\delta_{n,k}$, where $\delta$ denote the Kronecker symbol.
    
\end{lemme}

We can now define the renormalization of the $n$-th power of white noise images as a projection on the $n$-th chaos. Denote by $\Pi_n$ the orthogonal projection from $\mathbb{L}^2(\Omega,\mathcal{G})$ to $\mathcal{H}_n$.

\begin{lemme}

    Let $\rho\geqslant 0$ and $\eta\in\mathbb{L}^2(\mathbb{R}^2)$ such that $|\eta|_{\mathbb{L}^2_{x}}=1$ and define $Z=\rho W_\eta$. For all $n\in\N$, we have $\Pi_n(Z^n)=\rho^n\sqrt{n!}H_n(W_\eta).$
    
\end{lemme}

We define $\wick{Z^n}$ as $\Pi_n(Z^n)$. For example, in the case $n=2$ and $Z$ as in the previous lemma, it gives $\wick{Z^2}=Z^2-\rho^2=Z^2-\E\left[Z^2\right]$.\\ 

To prove the convergence of Wick product, first recall Nelson's estimate.

\begin{prop}{(Nelson's estimate)}\label{nelson-estimate}
    Let $n\in\N^*$ and $F\in\mathcal{H}_n$. Then, for all $p\geqslant 2$, we have $ \E\left[|F|^p\right] \leqslant (p-1)^{\frac{np}{2}} \E\left[|F|^{2}\right]^{\frac{p}{2}}$.

\end{prop}

The proof of \cref{conv_renorm} is very similar to the proof of Proposition 4 in \cite{debouard_2d_GP}. Thus, we only give here the technical argument specific to our case.

\begin{demo}{of \cref{conv_renorm}.}

    Let $i\in\{1,2\}$ and $N\in\N$, define $Z_{i,N} = \p_i Y_N$ where $Y_N$ was defined in \cref{subsection_noise_regularity}, 
    $$\forall x\in\mathbb{R}^2,\; \rho_N^i(x)^2 = \sum_{k\in\N^2} \chi_{k,N}^2 \frac{(\p_ih_k(x))^2}{\lambda_k^4} = \E\left[Z_{i,N}^2\right] \text{ and }\eta_N^i(x) = \frac{1}{\rho_N^i(x)}\sum_{k\in\N^2} \chi_{k,N} \frac{\p_ih_k(x)}{\lambda_k^2} h_k$$
    so that $|\eta_N^i(x)|_{\mathbb{L}^2_{\mathbb{R}^2}}=1$ and $Z_{i,N}(x)=\rho_N^i(x) W_{\eta_N^i(x)}$. Then, we can define $\wick{Z_{i,N}^2}=Z_{i,N}^2-\E\left[Z_{i,N}^2\right]$ so that $\wick{|\nabla Y_N|^2}=\wick{Z_{1,N}^2}+\wick{Z_{2,N}^2}$.\\ 

    We will prove the convergence of each $\wick{Z_{i,N}^2}$ separately. So fix $i\in\{1,2\}$ and let $q>2$ and $s>0$ which will be fixed later. As in the proof of proposition 4 in \cite{debouard_2d_GP}, we can show
    $$\E\left[\left|\wick{Z_{i,N}^2}-\wick{Z_{i,M}^2}\right|_{\W^{-s,q}_x}^q\right]\lesssim \left|K_{i,N,N}^2-2K_{i,N,M}^2+K_{i,M,M}^2\right|_{\mathbb{L}^{r}_{x_1}\W^{1-s+\delta,2}_{x_2}}^\frac{q}{2},$$
    with $0<\delta<s$, $\frac{2}{q}-\frac{s}{2}<\frac{1}{r}<\frac{2}{q}$ and for all $x,y\in\mathbb{R}^2$:
    $$K_{i,N,M}(x,y) = \sum_{k\in\N^2} \chi_{k,N}\chi_{k,M}\frac{\p_ih_k(x)\p_ih_k(y)}{\lambda_k^4}.$$
    Let $s'=s-\delta$. Using that for any $N,M\in\N^*$ we have 
    $$K_{i,N,N}^2-K_{i,N,M}^2 = (K_{i,N,N}-K_{i,N,M})(K_{i,N,N}+K_{i,N,M})$$ and applying \cref{product_rule_2var} with $r_1=r_2=2r$, $\alpha = 1 - s'$, $\sigma = 1-\frac{s'}{2}$ and $p=\frac{4}{s'}$ and the Sobolev embedding of $\W^{1-\frac{s'}{2},2}$ in $\mathbb{L}^{p}$ we obtain
    \begin{equation}\label{equation_conv_renorm}
        \begin{split}
            E\left[\left|\wick{Z_{i,N}^2}-\wick{Z_{i,M}^2}\right|_{\W^{-s,q}_x}^q\right]^{\frac{2}{q}} &\lesssim \left|K_{i,N,N}-K_{i,N,M}\right|_{\mathbb{L}^{2r}_{x}\W^{1-\frac{s'}{2},2}_{y}}\left|K_{i,N,N}+K_{i,N,M}\right|_{\mathbb{L}^{2r}_{x}\W^{1-\frac{s'}{2},2}_{y}}\\ 
            &+\left|K_{i,M,M}-K_{i,N,M}\right|_{\mathbb{L}^{2r}_{x}\W^{1-\frac{s'}{2},2}_{y}}\left|K_{i,M,M}+K_{i,N,M}\right|_{\mathbb{L}^{2r}_{x}\W^{1-\frac{s'}{2},2}_{y}}. 
        \end{split}
    \end{equation}
    Recall that Hermite functions verify the formula
    $$ \p_i h_k = \sqrt{\frac{k_i}{2}}h_{k-e_i}-\sqrt{\frac{k_i+1}{2}}h_{k+e_i}$$
    where $e_i=(\delta_{1,i},\delta_{2,i})$. This gives the following result:
    \begin{align*}
        K_{i,N,N}(x,y)-K_{i,N,M}(x,y) &= \sum_{k\in\N^2} \left[\sqrt{\frac{k_i+1}{2}}\frac{\chi_{k+e_i,N}(\chi_{k+e_i,N}-\chi_{k+e_i,M})\p_ih_{k+e_i}(x)}{\lambda_{k+e_i}^4}\right.\\ 
        &-\left.\sqrt{\frac{k_i}{2}}\frac{\chi_{k-e_i,N}(\chi_{k-e_i,N}-\chi_{k-e_i,M})\p_ih_{k-e_i}(x)}{\lambda_{k-e_i}^4}\right]h_k(y).
    \end{align*}
    Let $\alpha\in\R$, we deduce 
    \begin{equation*}
        |K_{i,N,N}(x,\cdot)-K_{i,N,M}(x,\cdot)|^2_{\W^{\alpha,2}_{y}}\lesssim \sum_{k\in\N^2} \frac{\chi_{k,N}^2(\chi_{k,N}-\chi_{k,M})^2}{\lambda_k^8}(\p_i h_k(x))^2 \left(k_i \lambda_{k-e_i}^{2\alpha}+(k_i+1) \lambda_{k+e_i}^{2\alpha}\right).
    \end{equation*}
    As $\lambda_k^2=2|k|+2$, we obtain 
    $$|K_{i,N,N}(x,\cdot)-K_{i,N,M}(x,\cdot)|^2_{\W^{\alpha,2}_{y}} \lesssim \sum_{\lfloor \frac{M}{2}\rfloor \leqslant |k| \leqslant N} \lambda_{k}^{2\alpha-6}(\p_i h_k(x))^2\lesssim \sum_{\lfloor \frac{M}{2}\rfloor-1 \leqslant |k| <+\infty} \lambda_{k}^{2\alpha-4}h_k(x)^2.$$
    Hence, for $M\geqslant 2$, we have
    $$|K_{i,N,N}(x,\cdot)-K_{i,N,M}(x,\cdot)|^2_{\W^{\alpha,2}_{y}} \lesssim \lambda_M^{-2\delta} \sum_{k\in\N^2} \lambda_{k}^{2\delta - 2\alpha-4}h_k(x)^2 \lesssim \lambda_M^{-2\delta} |K(x,\cdot)|_{\W^{\alpha+\delta,2}_{y}}^2.$$
    The same holds for $|K_{i,N,M}(x,\cdot)-K_{i,M,M}(x,\cdot)|^2_{\W^{\alpha,2}_{y}}$ and an analogous proof shows 
    $$\max\left(\left|K_{i,N,N}+K_{i,N,M}\right|_{\mathbb{L}^{2r}_{x}\W^{1-\frac{s'}{2},2}_{y}}, \left|K_{i,N,M}+K_{i,M,M}\right|_{\mathbb{L}^{2r}_{x}\W^{1-\frac{s'}{2},2}_{y}}\right) \lesssim |K|_{\mathbb{L}^{2r}_{x}\W^{1-\frac{s'}{2}+\delta,2}_{y}}.$$
    Going back to \cref{equation_conv_renorm}, we obtain
    $$\left|\wick{Z_{i,N}^2}-\wick{Z_{i,M}^2}\right|_{\mathbb{L}^q_\omega\W^{-s,q}_x}\lesssim \lambda_M^{-\delta}|K|_{\mathbb{L}^{2r}_{x}\W^{1-\frac{s'}{2}+\delta,2}_{y}}^2.$$
    Moreover, $|K|_{\mathbb{L}^{2r}_{x}\W^{1-\frac{s'}{2}+\delta,2}_{y}}$ is finite as soon as $1-\frac{s'}{2}+\delta < 1-\frac{1}{r}$ according to \cref{regu_kernel_-H-1}. There exists such a $\delta$ if $sq>2$. In order to have $0<s<1$, we need $q>2$. This proves the sequence is Cauchy. Then let $N$ go to infinity and the estimate follows. We can choose $\delta<\frac{s}{3}-\frac{2}{3r}$. Taking $\frac{1}{r}$ arbitrarily close to $\frac{2}{q}-\frac{s}{2}$, we obtain the announced range for $\delta$. \endproof
    
\end{demo}

\begin{rem}
    An analogous proof shows $\wick{|\nabla Y|^2}$ does not depend on the choice of $\chi\in\mathcal{S}(\R)$ (as long as $\chi(0)=1$), but choosing a $\chi$ which is compactly supported and equal to one on a neighborhood of the origin makes the proof of \cref{conv_renorm} easier and gives an explicit convergence speed. This independence, together with \cref{Proposition-unicity-nonlinear-equation} below, ensures that the renormalized solution we construct does not depend on the choice of $\chi$.
\end{rem}

To estimate $\mathbb{L}^q_x$ norms of Wick products, we use the following lemmas. First see that Nelson's estimate allows us to estimate higher order moments of $\W^{\alpha,q}_x$ norms of the Wick product.

\begin{lemme}\label{nelson-estimate-norm-W-alpha-q}
    Let $2\leqslant q<q' <+\infty$ and $\alpha\in\R$. It holds,
    $$\forall N\in\N,\; \left|\wick{|\nabla Y_N|^2}\right|_{\mathbb{L}^{q'}_\omega\W^{\alpha,q}_x} \leqslant (q'-1) \left|\wick{|\nabla Y_N|^2}\right|_{\mathbb{L}^{q}_\omega\W^{\alpha,q}_x}.$$
\end{lemme}

Using Borel-Cantelli lemma, we obtain the following corollary.

\begin{cor}\label{Cor:cauchy-criter-W(-kappa;infty)}
    Let $q\in(2,+\infty]$, $\frac{2}{q}<s<1$ and $0<\delta<\frac{2s}{3}-\frac{4}{3q}$. Almost surely $$\displaystyle\sup_{N\in\N} \left[\lambda_N^\delta\left|\wick{|\nabla Y_N|^2}-\wick{|\nabla Y|^2}\right|_{\W^{-s,q}_x} \right]<+\infty.$$
    
    Moreover, for any $p\in[1,+\infty)$ there exists a constant $C>0$ such that
    $$\forall N\geqslant M,\; \left|\wick{|\nabla Y_N|^2}-\wick{|\nabla Y_M|^2}\right|_{\mathbb{L}^p_\omega\W^{-s,q}_x} \leqslant C \lambda_M^{-\delta}$$
    and
    $$\max\left(\left|\wick{|\nabla Y|^2}\right|_{\mathbb{L}^p_\omega\W^{-s,q}_x},\sup_{N\in\N}\left|\wick{|\nabla Y_N|^2}\right|_{\mathbb{L}^p_\omega\W^{-s,q}_x}\right)\leqslant C.$$
\end{cor}

\begin{rem}
    This shows that, up to restricting ourselves to an event of probability 1, we can assume in what follows that for all $\omega\in\Omega$, $\wick{|\nabla Y_N|^2}(\omega)$ and $\wick{|\nabla Y|^2}(\omega)$ are well-defined and for every $q\in(2,+\infty]$ and $s>\frac{2}{q}$, $\wick{|\nabla Y_N|^2}(\omega)$ converges to $\wick{|\nabla Y|^2}(\omega)$ in $\W^{-s,q}$.\\ 
\end{rem}

Finally, in order to control the divergence of spatial $\mathbb{L}^q$-norms of Wick products, we will need the following lemma. 

\begin{lemme}\label{control-norm-Lq-Wick-product}

    Let $2\leqslant q <+\infty$ and $2\leqslant p_1,q_1,p_2,q_2<+\infty$ such that $\frac{1}{p_1}+\frac{1}{q_1}=\frac{1}{p_2}+\frac{1}{q_2}=\frac{1}{q}$. There exists $C>0$ such that
    $$\forall N\in\N,\;\left|\wick{|\nabla Y_N|^2}\right|_{\mathbb{L}^q_\omega\mathbb{L}^{q}_{x}} \leqslant C |\nabla Y_N|_{\mathbb{L}^{p_1}_\omega\mathbb{L}^{p_2}_{x}}|\nabla Y_N|_{\mathbb{L}^{q_1}_\omega\mathbb{L}^{q_2}_{x}}.$$
    
\end{lemme}

We deduce the following corollary.

\begin{cor}\label{estim-Lq-norm-Wick-product}

    Let $1\leqslant p <+\infty$, $2<q<+\infty$ and $1>s>\frac{2}{q}$. There exists $C>0$ such that for all $N\in\N$, we have
    $$\left|\wick{|\nabla Y_N|^2}\right|_{\mathbb{L}^p_\omega\mathbb{L}^q_{x}} \leqslant C \lambda_N^{s}.$$
    
\end{cor}

Hölder's inequality, Sobolev embeddings and \cref{Lemma_bound_from_moment} then give an immediate corollary.

\begin{cor}\label{estim-moment-norm-Lq-X.produit-wick}

    Let $\kappa\in\left(0,1\right)$, $q\in(1,+\infty)$ and a random variable $X\in\W^{1,q}$ a.s. There exists a sequence of random variables $(C_N)_{N\in\N}$ verifying \hyperlink{property_star}{Property (*)} such that almost surely 
    $$\forall N\in\N,\; |X \wick{|\nabla Y_N|^2}|_{\mathbb{L}^{q}_{x}}\leqslant C_N|X|_{\W^{1,q}_x} \lambda_N^{\kappa}.$$
    Moreover, let $p_0>p\geqslant 1$ and assume $X\in\mathbb{L}^{p_0}(\Omega,\W^{1,q})$, then there exists $C>0$ such that
    $$\sup_{N\in\N}\left[ \lambda_N^{-\kappa}|X \wick{|\nabla Y_N|^2}|_{\mathbb{L}^p_\omega\mathbb{L}^{q}_{x}}\right]\leqslant C |X|_{\mathbb{L}^{p_0}_\omega\W^{1,q}_x}.$$

\end{cor}

\subsection{The regularized equation}\label{subsection_regu_eq}

In order to solve \cref{Eq:2D-transformed-renorm} we are going to study a regularized transformed equation 
\begin{equation}{~}
    \begin{cases}\mathrm{i}\p_t v^N+ Hv^N +2\nabla Y_N\cdot \nabla v^N+ xY_N\cdot x v^N + \wick{|\nabla Y_N|^2} v^N +\lambda \left|v^N\right|^{2}v^N\e^{2 Y_N}=0\\ v(0)=v_0\end{cases} \label{Eq:2D-transformed-renorm-regularized}
\end{equation}
where $N\in\N$ and $Y_N = S_N Y$ as previously. Formally, this renormalized transformed equation has two invariant quantities, first a transformed mass 
\begin{equation}
    \Tilde{M}_N(v)=\displaystyle\int_{\mathbb{R}^2} |v(x)|^2 \e^{2Y_N(x)} \d x \label{Eq:def-trans-mass-regu}
\end{equation}
and a transformed energy 
\begin{equation}
    \Tilde{\mathcal{E}}_N(v)=\frac{1}{2}\displaystyle\int_{\mathbb{R}^2}\left( \left|\nabla v\right|^2+\left|xv\right|^2-\left|xv\right|^2Y_N-\wick{|\nabla Y_N|^2} \left|v\right|^2 \right)\e^{2Y_N}\d x  - \frac{\lambda}{4}\int_{\mathbb{R}^2}\left|v\right|^{4} \e^{4 Y_N}\d x.\label{Eq:def-trans-energy-regu}
\end{equation}

It is well-known that \cref{Eq:2D-transformed-renorm-regularized} has a unique maximal solution for $v_0\in\Sigma$ (see for example \cite{cazenave2003semilinear}). Here we give a complete statement. 

\begin{theo}{(Local wellposedness)}\label{Theorem-local-wellposedness}
    Let $\lambda\in\R$ and $N\in\N$. Almost surely, for every $v_0\in\Sigma$, there exists a unique maximal solution $v^N$ to \cref{Eq:2D-transformed-renorm-regularized} in $\mathcal{C}((-T_N^-,T_N^+),\Sigma)$.
    Moreover we have almost surely, the blow-up alternative
    $$\lim_{t\to T_N^{\pm}} |v^N(t)|_\Sigma = +\infty \text{ or } T_N^{\pm}=+\infty,$$
    transformed mass and energy, defined by \cref{Eq:def-trans-mass-regu,Eq:def-trans-energy-regu}, are conserved and if $v_0\in\W^{2,2}$, then $v^N\in\mathcal{C}((-T^N_-,T^N_+),\W^{2,2})\cap\mathcal{C}^1((-T^N_-,T^N_+),\mathbb{L}^2(\mathbb{R}^2))$.

\end{theo}

The mass conservation allows us to control the $\mathbb{L}^2$ norm of $v^N$:
$$|v^N(t)|_{\mathbb{L}^2_x}^2 \leqslant |\e^{-2Y_N}|_{\mathbb{L}^\infty_{x}}\Tilde{M}_N(v^N(t))\leqslant |\e^{-2Y_N}|_{\mathbb{L}^\infty_{x}}\Tilde{M}_N(v_0) \leqslant |\e^{-2Y_N}|_{\mathbb{L}^\infty_{x}}|\e^{2Y_N}|_{\mathbb{L}^\infty_{x}}|v_0|_{\mathbb{L}^2_{x}}.$$
Moreover the energy conservation allows us to control the $\Sigma$ norm of $v^N$ in $\mathbb{L}^p(\Omega)$ in the linear case  (for some $p\geqslant 2$ depending on the integrability of $v_0$), almost surely in the defocusing case and on a non-zero probability event in the focusing case.

\begin{prop}\label{Proposition_energy_bound_linear}

    Assume that $v_0\in\Sigma$ a.s. and $\lambda=0$. Then, for any $N\in\N$, $v^N$ is global almost surely. Moreover, let $p_0\in(2,+\infty)$ and assume $v_0\in\mathbb{L}^{p_0}(\Omega,\Sigma)$, then for any $p\in[2,p_0)$, there exists $C>0$ such that
    $$\forall N\in\N,\; \left|v^N\right|_{\mathbb{L}^p_\omega\mathbb{L}^\infty_t\Sigma}\leqslant C |v_0|_{\mathbb{L}^{p_0}_\omega\Sigma}.$$
    
\end{prop}

In the nonlinear case, we will need the following Gagliardo-Nirenberg inequality.

\begin{lemme}{(See \cite{BREZIS1980677})}\label{Lemma_Gagliardo_Nirenberg}
    For all $u\in W^{1,2}$, $|u|_{\mathbb{L}^4_{x}}^4\leqslant \frac{1}{2} |u|^2_{\mathbb{L}^2_{x}}|\nabla u|^2_{\mathbb{L}^2_{x}}.$
\end{lemme}

\begin{prop}\label{Proposition_energy_bound_defocusing}

    Assume $\lambda\leqslant0$. Then, there exists a positive random variable $C$ almost surely finite such that almost surely, for any $v_0\in\Sigma$ and any $N\in\N$, $v^N$ is global and 
    $$\forall N\in\N,\; \left|v^N\right|_{\mathbb{L}^\infty_t\Sigma}\leqslant C (1+|v_0|_{\Sigma})^2.$$
    
\end{prop}

\begin{prop}\label{Proposition_energy_bound_focusing}

    Assume $\lambda>0$ and let $L>0$. Let $\overline{\Omega}_L=\left\{\lambda L^2 \left|\e^{-2 Y}\right|_{\mathbb{L}^\infty_{x}}^2\left|\e^{2 Y}\right|_{\mathbb{L}^\infty_{x}}\left|\e^{4 Y}\right|_{\mathbb{L}^\infty_{x}}<4\right\}.$ Then on $\overline{\Omega}_L$ (endowed with the conditional probability and $\sigma$-algebra), there exists a random variable $N_0\in\N$ a.s. and a positive random variable $C$ almost surely finite such that almost surely, for any $v_0\in\Sigma$ with $|v_0|_{\mathbb{L}^2_{x}}\leqslant L$ and any $N\geqslant N_0$, $v^N$ is global and 
    $$\forall N\geqslant N_0,\; \left|v^N\right|_{\mathbb{L}^\infty_t\Sigma}\leqslant C|v_0|_\Sigma.$$
    
\end{prop}

Unfortunately, this bound is not sufficient to obtain a uniform bound on $(\p_t v^N)$ because of the negative regularity of the noisy terms.

\section{The linear case}\label{Section_linear_case}

In what follows, we denote by $v^N$ the unique global solution of the linear equation
\begin{equation}{~}
    \begin{cases}\mathrm{i}\p_t v^N+ Hv^N +2\nabla Y_N\cdot \nabla v^N+ xY_N\cdot x v^N + \wick{|\nabla Y_N|^2} v^N =0\\ v(0)=v_0\in\mathbb{L}^{p_0}\left(\Omega,\W^{2,2}\right)\end{cases} \label{Eq:2D-transformed-renorm-regularized-linear}
\end{equation}
with $p_0\in(2,+\infty)$. First, as in \cite{debussche_weber_T2}, we prove a diverging bound on $v^N$ in $\W^{2,2}$.

\begin{prop}\label{Prop:W22-bound}

    Let $\kappa\in(0,1)$ and $2\leqslant p<p_0$. There exists a constant $C>0$ such that for any $N\in\N$,
    $$|v^N|_{\mathbb{L}^p_\omega\mathbb{L}^\infty_t\W^{2,2}_x} \leqslant C  \lambda_N^\kappa|v_0|_{\mathbb{L}^{p_0}_\omega\W^{2,2}_x}$$
    
\end{prop}

\begin{demo}
    
Define $w^N=\p_t v^N$. Using \cref{Eq:2D-transformed-renorm-regularized-linear}, we have 
$$Hv^N = -\left(\mathrm{i} w^N + 2\nabla Y_N\cdot \nabla v^N+ xY_N\cdot x v^N + \wick{|\nabla Y_N|^2} v^N\right).$$
Taking the $\mathbb{L}^2_x$ norm, we obtain
$$|v^N|_{\W^{2,2}_x} \leqslant \left|w^N\right|_{\mathbb{L}^2_{x}}+2\left|\nabla Y_N\cdot \nabla v^N\right|_{\mathbb{L}^2_{x}}+\left|x Y_N\cdot x v^N\right|_{\mathbb{L}^2_{x}}+\left|\wick{|\nabla Y_N|^2} v^N \right|_{\mathbb{L}^2_{x}}.$$
Applying \cref{Cor:estimate-norm-Lq-X-Y_N,estim-moment-norm-Lq-X.produit-wick,Proposition_energy_bound_linear}, we have for any $p<p'<p_0$
\begin{equation}
    |v^N|_{\mathbb{L}^p_\omega\mathbb{L}^\infty_t\W^{2,2}_x} \lesssim \left|w^N\right|_{\mathbb{L}^p_\omega\mathbb{L}^\infty_t\mathbb{L}^2_{x}} + \lambda_N^\kappa|v^N|_{\mathbb{L}^{p'}_\omega\mathbb{L}^\infty_t\Sigma} \lesssim \left|w^N\right|_{\mathbb{L}^p_\omega\mathbb{L}^\infty_t\mathbb{L}^2_{x}} + \lambda_N^\kappa|v_0|_{\mathbb{L}^{p_0}_\omega\Sigma}.\label{Equation_diverging_bound_linear_wN+v0Sigma}
\end{equation}
To bound the $\mathbb{L}^2_x$ norm of $w^N$, we use the fact that $w^N$ verifies at least formally \cref{Eq:2D-transformed-renorm-regularized-linear} with initial data $\mathrm{i}\left(Hv_0+2\nabla v_0\cdot\nabla Y_N+xv_0\cdot xY_N+\wick{|\nabla Y_N|^2}v_0\right)\in\mathbb{L}^2(\mathbb{R}^2)$ a.s. Approximating $v_0$ by a sequence of Schwartz functions, one can show that $w^N$ has its transformed mass conserved
$$\forall t\in\R,\; \int_{\mathbb{R}^2} |w^N(t)|\e^{2Y_N}\d x=\int_{\mathbb{R}^2} |w^N(0)|\e^{2Y_N}\d x.$$
Thus, by Hölder's inequality and \cref{Cor:moment-boud-expYN-norm-Linf}, we have for any $p<p'<p_0$:
\begin{align*}
    |w^N|_{\mathbb{L}^{p}_\omega\mathbb{L}^\infty_t\mathbb{L}^2_{x}} &\lesssim|w^N(0)|_{\mathbb{L}^{p'}_\omega\mathbb{L}^2_{x}}\\ 
    &\lesssim|Hv_0|_{\mathbb{L}^{p'}_\omega\mathbb{L}^2_{x}}+2|\nabla Y_N\cdot\nabla v_0|_{\mathbb{L}^{p'}_\omega\mathbb{L}^2_{x}}+|x Y_N\cdot xv_0|_{\mathbb{L}^{p'}_\omega\mathbb{L}^2_{x}}+|\wick{|\nabla Y_N|^2}v_0|_{\mathbb{L}^{p'}_\omega\mathbb{L}^2_{x}}\\ 
    &\lesssim |v_0|_{\mathbb{L}^{p_0}_\omega\W^{2,2}_x}+\lambda_N^\kappa|v_0|_{\mathbb{L}^{p_0}_\omega\Sigma}
\end{align*}
using \cref{Cor:estimate-norm-Lq-X-Y_N,estim-moment-norm-Lq-X.produit-wick}. Going back to \cref{Equation_diverging_bound_linear_wN+v0Sigma}, we finally obtain 
$$|v^N|_{\mathbb{L}^p_\omega\mathbb{L}^\infty_t\W^{2,2}_x} \lesssim  |v_0|_{\mathbb{L}^{p_0}_\omega\W^{2,2}_x} + \lambda_N^\kappa|v_0|_{\mathbb{L}^{p_0}_\omega\Sigma} \lesssim \lambda_N^\kappa|v_0|_{\mathbb{L}^{p_0}_\omega\W^{2,2}_x}.$$\endproof
\end{demo}

In order to manage this divergence, we prove a Cauchy estimate in $\mathbb{L}^p(\Omega,\mathcal{C}([a,b],\mathbb{L}^2(\mathbb{R}^2)))$. 

\begin{prop}\label{norm-L2-diff-NM}

    Let $N\geqslant M>0$ be two integers, $p\in[2,p_0)$, $\delta\in(0,\frac{2}{3})$ and $a<0<b$. Then there exists a constant $C>0$ such that we have
    $$|v^N-v^M|_{\mathbb{L}^p_\omega\mathbb{L}^\infty_{t}\mathbb{L}^2_{x}}^2 \leqslant C \lambda_M^{-\delta}\left|v_0\right|_{\mathbb{L}^{p_0}_\omega\W^{2,2}_x}^2.$$

\end{prop}

\begin{demo}
    
Let $N\geqslant M>0$ and set $R = v^N-v^M$. It verifies
\begin{equation}{~}
    \begin{cases}-\mathrm{i}\p_t R =& HR+2\nabla Y_N\cdot \nabla R+2\nabla (Y_N-Y_M)\cdot \nabla v^M+ xY_N\cdot x R\\ & + x(Y_N-Y_M)\cdot x v^M + \wick{|\nabla Y_N|^2} R+ \left(\wick{|\nabla Y_N|^2}-\wick{|\nabla Y_M|^2} \right)v^M\\  R(0)=0\end{cases} \label{2d_transf_renorm_regu_eq-lin-diff-N-M}
\end{equation}
and it belongs to $\mathcal{C}^1([a,b],\mathbb{L}^2(\mathbb{R}^2))$ almost surely. Taking the time derivative of $\Tilde{M}_N(R)$, we obtain
\begin{align*}
    \frac{1}{2}\frac{\d}{\d t}\int_{\mathbb{R}^2} |R|^2\e^{2Y_N}\d x &= \langle \p_t R,R\e^{2Y_N}\rangle_{\mathbb{L}^2}\\ 
    &= - \im \int_{\mathbb{R}^2} \left(HR+2\nabla Y_N\cdot \nabla R+2\nabla (Y_N-Y_M)\cdot \nabla v^M \right)\overline{R}\e^{2Y_N}\d x\\ 
    &- \im \int_{\mathbb{R}^2} \left(xY_N\cdot x R + x(Y_N-Y_M)\cdot x v^M\right)\overline{R}\e^{2Y_N}\d x\\ 
    & - \im \int_{\mathbb{R}^2} \left(\wick{|\nabla Y_N|^2} R+ \left(\wick{|\nabla Y_N|^2}-\wick{|\nabla Y_M|^2} \right)v^M\right)\overline{R}\e^{2Y_N}\d x\\ 
    &=- \im \int_{\mathbb{R}^2} \left(2\nabla (Y_N-Y_M)\cdot \nabla v^M+ x(Y_N-Y_M)\cdot x v^M\right)\overline{R}\e^{2Y_N}\d x\\ 
    &-\im \int_{\mathbb{R}^2}\left(\wick{|\nabla Y_N|^2}-\wick{|\nabla Y_M|^2} \right)v^M\overline{R}\e^{2Y_N}\d x,
\end{align*}
where we used 
\begin{align*}
    \im \int_{\mathbb{R}^2} \left(HR+2\nabla Y_N\cdot \nabla R\right)\overline{R}\e^{2Y_N}\d x &= \im \int_{\mathbb{R}^2}\left( \overline{R}\e^{2Y_N} \Delta R + \overline{R} \nabla R\cdot\nabla\left[\e^{2Y_N}\right]\right) \d x\\ 
    &= \im \int_{\mathbb{R}^2}\left( -\nabla R\cdot \nabla\left[\overline{R}\e^{2Y_N}\right] + \overline{R} \nabla R\cdot\nabla\left[\e^{2Y_N}\right]\right) \d x\\ 
    &=0.
\end{align*}
Let $\kappa\in(0,1)$ such that $\delta\in(0,\frac{2\kappa}{3})$. We have by duality:
\begin{align*}
    \frac{1}{2}\frac{\d}{\d t}\int_{\mathbb{R}^2} |R|^2\e^{2Y_N} &\leqslant 2 |\nabla(Y_N-Y_M)|_{\W^{-\kappa,\infty}_x}\left|\overline{R}\nabla v^M\e^{2Y_N}\right|_{\W^{\kappa,1}_x} \\ 
    &+ |x(Y_N-Y_M)|_{\W^{-\kappa,\infty}_x}\left|x \overline{R} v^M\e^{2Y_N}\right|_{\W^{\kappa,1}_x} \\ 
    &+ \left|\wick{|\nabla Y_N|^2}-\wick{|\nabla Y_M|^2}\right|_{\W^{-\kappa,\infty}_x}\left|\overline{R} v^M\e^{2Y_N}\right|_{\W^{\kappa,1}_x}.
\end{align*}
By \cref{product_rule} and interpolation, we have
\begin{align*}
    \left|\overline{R}\nabla v^M\e^{2Y_N}\right|_{\W^{\kappa,1}_x} &\lesssim \left|R\right|_{\W^{\kappa,2}_x}\left|v^M\right|_{\W^{1+\kappa,2}_x}\left|\e^{2Y_N}\right|_{1+\W^{\kappa,\infty}_x}\\ 
    &\lesssim \left(|v^N|_{\mathbb{L}^\infty_t\Sigma} +|v^M|_{\mathbb{L}^\infty_t\Sigma}\right)|v^M|_{\mathbb{L}^\infty_t\Sigma}^{1-\kappa}\left|\e^{2Y_N}\right|_{1+\W^{\kappa,\infty}_x}\left|v^M\right|_{\mathbb{L}^\infty_t\W^{2,2}_x}^\kappa
\end{align*}
and
\begin{align*}
    |x\overline{R} v^M\e^{2Y_N}|_{\W^{\kappa,1}_x} &\lesssim \left(|v^N|_{\mathbb{L}^\infty_t\Sigma} +|v^M|_{\mathbb{L}^\infty_t\Sigma}\right)|v^M|_{\mathbb{L}^\infty_t\Sigma}^{1-\kappa}\left|\e^{2Y_N}\right|_{1+\W^{\kappa,\infty}_x}\left|v^M\right|_{\mathbb{L}^\infty_t\W^{2,2}_x}^\kappa
\end{align*}
and
\begin{align*}
    |\overline{R}v^M\e^{2Y_N}|_{\W^{\kappa,1}_x}&\lesssim \left(|v^N|_{\mathbb{L}^\infty_t\Sigma} +|v^M|_{\mathbb{L}^\infty_t\Sigma}\right)|v^M|_{\mathbb{L}^\infty_t\Sigma}\left|\e^{2Y_N}\right|_{1+\W^{\kappa,\infty}_x}.
\end{align*}
Integrating in time, we obtain a bound uniform in $t\in[a,b]$ (with a constant depending on $b-a$):
\begin{align*}
    |R|_{\mathbb{L}^\infty_{t}\mathbb{L}^2_{x}}^2 &\lesssim \left|e^{-2Y_N}\right|_{\mathbb{L}^\infty_{x}}\left[|\nabla(Y_N-Y_M)|_{\W^{-\kappa,\infty}_x}|v^M|_{\mathbb{L}^\infty_t\Sigma}^{1-\kappa}\left|\e^{2Y_N}\right|_{1+\W^{\kappa,\infty}_x}\left|v^M\right|_{\mathbb{L}^\infty_t\W^{2,2}_x}^\kappa\right. \\ 
    &+ |x(Y_N-Y_M)|_{\W^{-\kappa,\infty}_x}|v^M|_{\mathbb{L}^\infty_t\Sigma}^{1-\kappa}\left|\e^{2Y_N}\right|_{1+\W^{\kappa,\infty}_x}\left|v^M\right|_{\mathbb{L}^\infty_t\W^{2,2}_x}^\kappa\\ 
    &+ \left.\left|\wick{|\nabla Y_N|^2}-\wick{|\nabla Y_M|^2}\right|_{\W^{-\kappa,\infty}_x}|v^M|_{\mathbb{L}^\infty_t\Sigma}\left|\e^{2Y_N}\right|_{1+\W^{\kappa,\infty}_x}\right]\left(|v^N|_{\mathbb{L}^\infty_t\Sigma} +|v^M|_{\mathbb{L}^\infty_t\Sigma}\right).
\end{align*}
Let $p<p'<p_0$ and $r>1$ such that $\frac{2}{p}>\frac{2}{p'}+\frac{1}{r}$, by Hölder's inequality and \cref{Cor:moment-boud-expYN-norm-Linf} we have 
\begin{align*}
    |R|_{\mathbb{L}^p_\omega\mathbb{L}^\infty_{t}\mathbb{L}^2_{x}}^2 &\lesssim |\nabla(Y_N-Y_M)|_{\mathbb{L}^{r}_\omega\W^{-\kappa,\infty}_x}\left(|v^N|_{\mathbb{L}^{p'}_\omega\mathbb{L}^\infty_t\Sigma} +|v^M|_{\mathbb{L}^{p'}_\omega\mathbb{L}^\infty_t\Sigma}\right)|v^M|_{\mathbb{L}^{p'}_\omega\mathbb{L}^\infty_t\Sigma}^{1-\kappa}\left|v^M\right|_{\mathbb{L}^{p'}_\omega\mathbb{L}^\infty_t\W^{2,2}_x}^\kappa \\ 
    &+ |x(Y_N-Y_M)|_{\mathbb{L}^{r}_\omega\W^{-\kappa,\infty}_x}\left(|v^N|_{\mathbb{L}^{p'}_\omega\mathbb{L}^\infty_t\Sigma} +|v^M|_{\mathbb{L}^{p'}_\omega\mathbb{L}^\infty_t\Sigma}\right)|v^M|_{\mathbb{L}^{p'}_\omega\mathbb{L}^\infty_t\Sigma}^{1-\kappa}\left|v^M\right|_{\mathbb{L}^{p'}_\omega\mathbb{L}^\infty_t\W^{2,2}_x}^\kappa\\ 
    &+ \left|\wick{|\nabla Y_N|^2}-\wick{|\nabla Y_M|^2}\right|_{\mathbb{L}^{r}_\omega\W^{-\kappa,\infty}_x}\left(|v^N|_{\mathbb{L}^{p'}_\omega\mathbb{L}^\infty_t\Sigma} +|v^M|_{\mathbb{L}^{p'}_\omega\mathbb{L}^\infty_t\Sigma}\right)|v^M|_{\mathbb{L}^{p'}_\omega\mathbb{L}^\infty_t\Sigma}.
\end{align*}
Let $\kappa'\in(0,1)$, using \cref{Proposition_energy_bound_linear,Prop:W22-bound}, we obtain 
\begin{align*}
    |R|_{\mathbb{L}^p_\omega\mathbb{L}^\infty_{t}\mathbb{L}^2_{x}}^2 &\lesssim \left(|\nabla(Y_N-Y_M)|_{\mathbb{L}^{r}_\omega\W^{-\kappa,\infty}_x}+|x(Y_N-Y_M)|_{\mathbb{L}^{r}_\omega\W^{-\kappa,\infty}_x}\right)|v_0|^{2-\kappa}_{\mathbb{L}^{p_0}_\omega\Sigma}|v_0|^{\kappa}_{\mathbb{L}^{p_0}_\omega\W^{2,2}}\lambda_M^{\kappa \kappa'}\\ 
    &+ \left|\wick{|\nabla Y_N|^2}-\wick{|\nabla Y_M|^2}\right|_{\mathbb{L}^{r}_\omega\W^{-\kappa,\infty}_x}|v_0|^{2}_{\mathbb{L}^{p_0}_\omega\Sigma}.
\end{align*}
Let $\kappa''\in(\delta,\kappa)$, using \cref{Cor:convergence-Y_N-W(1-kappa;infty)}, we obtain
\begin{align*}
    |R|_{\mathbb{L}^p_\omega\mathbb{L}^\infty_{t}\mathbb{L}^2_x}^2 &\lesssim |v_0|^{2-\kappa}_{\mathbb{L}^{p_0}_\omega\Sigma}|v_0|^{\kappa}_{\mathbb{L}^{p_0}_\omega\W^{2,2}_x}\lambda_M^{\kappa \kappa'-\kappa''}+ \left|\wick{|\nabla Y_N|^2}-\wick{|\nabla Y_M|^2}\right|_{\mathbb{L}^{r}_\omega\W^{-\kappa,\infty}_x}|v_0|^{2}_{\mathbb{L}^{p_0}_\omega\Sigma}.
\end{align*}
Using \cref{Cor:cauchy-criter-W(-kappa;infty)} (recall that $\delta\in(0,\frac{2\kappa}{3})$), we have
\begin{align*}
    |R|_{\mathbb{L}^p_\omega\mathbb{L}^\infty_{t}\mathbb{L}^2_{x}}^2 &\lesssim |v_0|_{\mathbb{L}^{p_0}_\omega\W^{2,2}_x}^2\left(\lambda_M^{\kappa \kappa'-\kappa''}+\lambda_M^{-\delta}\right).
\end{align*}
We conclude by choosing $\kappa' = \frac{\kappa''-\delta}{\kappa}\in(0,1)$.\endproof

\end{demo}

We can now prove our main result in the linear case.

\begin{demo}{of \cref{Theo-lin}.}

    Let $\kappa'\in(0,1)$ and $\delta\in(0,\frac{2}{3})$. By interpolation between \cref{Prop:W22-bound,norm-L2-diff-NM} we get
    $$|v^N-v^M|_{\mathbb{L}^p_\omega\mathbb{L}^\infty_{t}\W^{\sigma,2}_x} \lesssim |v^N-v^M|_{\mathbb{L}^p_\omega\mathbb{L}^\infty_{t}\W^{2,2}_x}^{\frac{\sigma}{2}}|v^N-v^M|_{\mathbb{L}^p_\omega\mathbb{L}^\infty_{t}\mathbb{L}^{2}_{x}}^{1-\frac{\sigma}{2}}\lesssim \lambda_N^{\frac{\kappa'\sigma}{2}}\lambda_M^{-\left(1-\frac{\sigma}{2}\right)\frac{\delta}{2}}\left|v_0\right|_{\mathbb{L}^{p_0}_\omega\W^{2,2}_x}.$$
    Let $n\in\N$, $N=2^{n+1}$ and $M=2^n$. We have
    $$\left|v^{2^{n+1}}-v^{2^n}\right|_{\mathbb{L}^p_\omega\mathbb{L}^\infty_{t}\W^{\sigma,2}_x} \lesssim 2^{n\left(\kappa'\sigma-\left(1-\frac{\sigma}{2}\right)\delta\right)/4}\left|v_0\right|_{\mathbb{L}^{p_0}_\omega\W^{2,2}_x}. $$
    Remark that:
    $$\frac{\kappa'\sigma}{4}-\left(1-\frac{\sigma}{2}\right)\frac{\delta}{4} < 0 \iff \kappa' < \frac{2-\sigma}{2\sigma}\delta.$$
    Thus, for fixed $\sigma$, it is possible to choose  $\delta\in(0,\frac{1}{4})$ and $\kappa'\in\left(0,\frac{2-\sigma}{2\sigma}\delta\right)$ so that the sequence $\left(\left|v^{2^{n+1}}-v^{2^n}\right|_{\mathbb{L}^p_\omega\mathbb{L}^\infty_{t}\W^{\sigma,2}_x}\right)_{n\in\N}$ is summable. Hence, $\left(v^{2^n}\right)_{n\in\N}$ is a Cauchy sequence in the Banach space $\mathbb{L}^p\left(\Omega,\mathcal{C}([a,b],\W^{\sigma,2})\right)$. Denote by $v$ its limit in this space. Using \cref{product-rule-W(-k;r)-Wkp,Lem:Conv-W(1-kappa;inf)-exp(aYN),Cor:cauchy-criter-W(-kappa;infty),Cor:convergence-Y_N-W(1-kappa;infty)}, it is not difficult to show that $-\mathrm{i}\p_t v^{2^n}$ converge in $\mathcal{C}([a,b],\W^{\sigma-2,2})$ to $$Hv+2\nabla v\cdot \nabla Y + xv\cdot xY+\wick{|\nabla Y|^2} v $$ so that $v\in\mathcal{C}^1([a,b],\W^{\sigma-2,2})\cap\mathcal{C}([a,b],\W^{\sigma,2})$ a.s. and is solution of \cref{Eq:2D-transformed-renorm-linear}.\\ 

    Since \cref{Eq:2D-transformed-renorm-linear} is linear, it is sufficient, in order to prove uniqueness, to prove that any solution with initial data $0$ is constant to $0$. Let $v$ be a solution in $\mathbb{L}^p\left(\Omega,\mathcal{C}([a,b],\W^{\sigma,2})\right)$ of \cref{Eq:2D-transformed-renorm-linear} with initial data 0. Then the Lions-Magenes lemma (see \cite{Lions1972-fp}) yields:
    \begin{align*}
        \frac{1}{2}\frac{\d}{\d t} \int_{\mathbb{R}^2} |v|^2 \e^{2Y_N}\d x &= \left\langle \p_t v, v\e^{2Y_N}\right\rangle_{\W^{\sigma-2,2}_x,\W^{2-\sigma,2}_x}\\ 
        &= \left\langle \mathrm{i}\left(Hv+2\nabla v\cdot\nabla Y+xv\cdot xY+:|\nabla Y|^2:v\right), v\e^{2Y_N}\right\rangle_{\W^{\sigma-2,2}_x,\W^{2-\sigma,2}_x}.
    \end{align*}
    Recall that
    $$\left\langle \mathrm{i}\left(Hv+2\nabla v\cdot\nabla Y_N+xv\cdot xY_N+:|\nabla Y_N|^2:v\right), v\e^{2Y_N}\right\rangle_{\W^{\sigma-2,2}_x,\W^{2-\sigma,2}_x}=0.$$
    It follows
    \begin{align*}
        \frac{1}{2}\frac{\d}{\d t} \int_{\mathbb{R}^2} |v|^2 \e^{2Y_N}\d x &= \left\langle \mathrm{i}\left(2\nabla v\cdot\nabla (Y-Y_N)+xv\cdot x(Y-Y_N)v\right), v\e^{2Y_N}\right\rangle_{\W^{\sigma-2,2}_x,\W^{2-\sigma,2}_x}\\ 
        &+\left\langle \mathrm{i}\left(:|\nabla Y|^2:-:|\nabla Y_N|^2:\right)v, v\e^{2Y_N}\right\rangle_{\W^{\sigma-2,2}_x,\W^{2-\sigma,2}_x}.
    \end{align*}
    Let $\kappa\in(0,\sigma-1)$, we have a.s.
    \begin{align*}
        \left|\frac{1}{2}\frac{\d}{\d t} \int_{\mathbb{R}^2} |v|^2 \e^{2Y_N}\d x\right| &\lesssim\left(\left|\nabla(Y-Y_N)\right|_{\W^{-\kappa,\infty}_x}\left|v\right|_{\W^{1+\kappa,2}_x}+\left|x(Y-Y_N)\right|_{\W^{-\kappa,\infty}_x}\left|v\right|_{\W^{1+\kappa,2}_x}\right.\\ 
        &+\left.\left|:|\nabla Y|^2:-:|\nabla Y_N|^2:\right|_{\W^{-\kappa,\infty}_x}\left|v\right|_{\W^{\kappa,2}_x}\right)\left|v\right|_{\W^{\kappa,2}_x}\left|\e^{2Y_N}\right|_{1+\W^{\kappa,\infty}_x}\\ 
        &\lesssim\left(\left|\nabla(Y-Y_N)\right|_{\W^{-\kappa,\infty}_x}+\left|x(Y-Y_N)\right|_{\W^{-\kappa,\infty}_x}\right)\left|v\right|_{\mathbb{L}^\infty_{t}\W^{\sigma,2}_x}^2\left|\e^{2Y_N}\right|_{1+\W^{\kappa,\infty}_x}\\ 
        &+\left|:|\nabla Y|^2:-:|\nabla Y_N|^2:\right|_{\W^{-\kappa,\infty}_x}\left|v\right|_{\mathbb{L}^\infty_{t}\W^{\sigma,2}_x}^2\left|\e^{2Y_N}\right|_{1+\W^{\kappa,\infty}_x}.
    \end{align*}
    Integrating in time, taking expectation and using \cref{Lem:Conv-W(1-kappa;inf)-exp(aYN),Cor:cauchy-criter-W(-kappa;infty),Cor:convergence-Y_N-W(1-kappa;infty)} we obtain
    $$\sup_{[a,b]}\int_{\mathbb{R}^2} |v|^2 \e^{2Y} = 0 \text{ a.s.}$$
    which concludes the uniqueness.\\ 

    To prove now that the whole sequence converges to $v$, let $N\geqslant 2$ and $n=\lfloor\log_2(N)\rfloor$. Note $M=2^n$. Then, interpolating between \cref{Prop:W22-bound,norm-L2-diff-NM}, we have
    \begin{align*}
        |v^N-v|_{\mathbb{L}^p_\omega\mathbb{L}^\infty_{t}\W^{\sigma,2}_x} &\lesssim |v^N-v^{2^n}|_{\mathbb{L}^p_\omega\mathbb{L}^\infty_{t}\W^{2,2}_x}^{\frac{\sigma}{2}}|v^N-v^{2^n}|_{\mathbb{L}^p_\omega\mathbb{L}^\infty_{t}\mathbb{L}^{2}_{x}}^{1-\frac{\sigma}{2}}+|v^{2^n}-v|_{\mathbb{L}^p_\omega\mathbb{L}^\infty_{t}\W^{\sigma,2}_x}\\ 
        &\lesssim \lambda_N^{\frac{\kappa'\sigma}{2}}\lambda_M^{-\left(1-\frac{\sigma}{2}\right)\frac{\delta}{2}}+|v^{2^n}-v|_{\mathbb{L}^p_\omega\mathbb{L}^\infty_{t}\W^{\sigma,2}_x}\\ 
        &\lesssim 2^{\log_2(N)\frac{\kappa'\sigma}{4}-n\left(1-\frac{\sigma}{2}\right)\frac{\delta}{4}}+|v^{2^n}-v|_{\mathbb{L}^p_\omega\mathbb{L}^\infty_{t}\W^{\sigma,2}_x},
    \end{align*}
    which goes to 0 when $N$ goes to $+\infty$ if, for example, we choose $\delta\in\left(0,\frac{1}{4}\right)$, $\delta'\in\left(0,\left(1-\frac{\sigma}{2}\right)\frac{\delta}{4}\right)$ and $\kappa'\in\left(0,\min\left(1,\frac{(2-\sigma)\delta-8\delta'}{4\sigma}\right)\right)$ in order to have
    $$\frac{\log_2(N)}{\lfloor \log_2(N)\rfloor}\frac{\kappa'\sigma}{4}-\left(1-\frac{\sigma}{2}\right)\frac{\delta}{4}<-\delta'$$ uniformly in $N\geqslant 2$.\\ 
    
    It remains to prove conservation laws. The conservation of the transformed mass follows from the convergence of $v^N$ in $\mathbb{L}^2(\mathbb{R}^2)$ and \cref{Lem:Conv-W(1-kappa;inf)-exp(aYN)} and the conservation of the transformed energy follows from convergence of $v^N$ in $\Sigma$ and \cref{product_rule,product-rule-W(-k;r)-Wkp,Lem:Conv-W(1-kappa;inf)-exp(aYN)}. This conclude the proof.\endproof
    
\end{demo}

\section{The nonlinear equation}\label{Section_nonlinear_case}

First we prove that there exists at most one solution to \cref{Eq:2D-transformed-renorm}.

\begin{prop}\label{Proposition-unicity-nonlinear-equation}

    Let $\lambda\in\R$. Almost surely, for every $a<0<b$, every $\sigma\in(1,2)$ and every $v_0\in\W^{\sigma,2}$, there exists at most one solution to \cref{Eq:2D-transformed-renorm} in $\mathcal{C}\left([a,b],\W^{\sigma,2}\right)$ and it belongs to $\mathcal{C}\left([a,b],\W^{\sigma,2}\right)\cap\mathcal{C}^1\left([a,b],\W^{\sigma-2,2}\right)$.
    
\end{prop}

\begin{demo}

    Using the continuous injection of $\mathcal{C}([a',b'],\W^{\sigma',2})$ in $\mathcal{C}([a,b],\W^{\sigma,2})$ for $a'<a$, $b'>b$ and $\sigma'>\sigma$, it is sufficient to prove the claim almost surely for fixed $a$, $b$ and $\sigma$. So let $a<0<b$ and $\sigma\in(1,2)$.\\ 

    Let $\kappa\in(0,\sigma-1)$ and $\delta\in\left(0,\frac{2\kappa}{3}\right)$. Define 
    \begin{align*}
        \Omega_0&=\left\{e^{-2Y}\in\mathbb{L}^\infty(\mathbb{R}^2), \e^{2 Y_N}\xrightarrow[N\to+\infty]{\mathbb{L}^\infty(\mathbb{R}^2)}\e^{2Y},\sup_{N\in\N}\left|e^{2Y_N}-1\right|_{\W^{1-\kappa,\infty}_x}<+\infty\right\}\\ 
        \cap&\left\{\sup_{N\in\N}\lambda_N^{\delta}\max\left(|\nabla (Y-Y_N)|_{\W^{-\kappa,\infty}_x},|x(Y-Y_N)|_{\W^{-\kappa,\infty}_x},|\wick{|\nabla Y|^2}-\wick{|\nabla Y_N|^2}|_{\W^{-\kappa,\infty}_x}\right)<+\infty\right\},
    \end{align*}
    using \cref{Lem:Conv-W(1-kappa;inf)-exp(aYN),Cor:convergence-Y_N-W(1-kappa;infty),Cor:cauchy-criter-W(-kappa;infty)}, it follows that $\P(\Omega_0)=1$.\\ 
    
    Let $v_0\in\W^{\sigma,2}$, $v$ and $w$ be two solutions of \cref{Eq:2D-transformed-renorm} in $\mathcal{C}([a,b],\W^{\sigma,2})$ with initial data $v_0$ and set $r=v-w$. Then $r\in\mathcal{C}([a,b],\W^{\sigma,2})$ and verifies 
    \begin{equation}{~}
    \begin{cases}\mathrm{i}\p_t r + Hr+2\nabla Y\cdot \nabla r+xY\cdot x r + \wick{|\nabla Y|^2} r+ \lambda\left(|v|^{2}v-|w|^{2}w \right)e^{2 Y}\\ r(0)=0.\end{cases} 
    \end{equation}
    Moreover, using \cref{product-rule-W(-k;r)-Wkp}, one can show that $v,w\in\mathcal{C}\left([a,b],\W^{\sigma,2}\right)\cap\mathcal{C}^1\left([a,b],\W^{\sigma-2,2}\right)$ so that $r\in\mathcal{C}\left([a,b],\W^{\sigma,2}\right)\cap\mathcal{C}^1\left([a,b],\W^{\sigma-2,2}\right)$. Using \cref{Lem:Conv-W(1-kappa;inf)-exp(aYN)}, on $\Omega_0$, we have 
    \begin{equation}
        \sup_{t\in[a,b]}|r(t)|_{\mathbb{L}^2_{x}}^2 \leqslant \left|e^{-2Y}\right|_{\mathbb{L}^\infty_{x}} \sup_{t\in[a,b]}\int_{\mathbb{R}^2} |r|^2\e^{2Y}\d x \leqslant \left|e^{-2Y}\right|_{\mathbb{L}^\infty_{x}} \lim_{N\to+\infty} \sup_{t\in[a,b]}\int_{\mathbb{R}^2} |r|^2\e^{2Y_N}\d x.\label{Equation_unicity_control}
    \end{equation}
    Recall we chose $\kappa\in(0,\sigma-1)$. As in the linear case, by taking the time derivative, the Lions-Magenes lemma yields 
    \begin{align*}
        \left|\frac{1}{2}\frac{\d}{\d t} \int_{\mathbb{R}^2} |r|^2 \e^{2Y_N}\d x\right|&= \left| \left\langle \mathrm{i}\left(2\nabla r\cdot\nabla (Y-Y_N)+xr\cdot x(Y-Y_N)+\left(\wick{|\nabla Y|^2}-\wick{|\nabla Y_N|^2}\right)r\right.\right.\right.\\ 
        &\left.\left.\left.+\lambda\left(|v|^{2}v-|w|^{2}w \right)e^{2 Y}\right), r\e^{2Y_N}\right\rangle_{\W^{-\kappa,2}_x,\W^{\kappa,2}_x}\right|\\ 
        &\leqslant \left|2\nabla r\cdot\nabla (Y-Y_N)+xr\cdot x(Y-Y_N)\right|_{\W^{-\kappa,2}_x}\left|r\e^{2Y_N}\right|_{\W^{\kappa,2}_x}\\ 
        &+ \left|\left(\wick{|\nabla Y|^2}-\wick{|\nabla Y_N|^2}\right)r\right|_{\W^{-\kappa,2}_x}\left|r\e^{2Y_N}\right|_{\W^{\kappa,2}_x}\\ 
        &+\left|\lambda\int_{\mathbb{R}^2}\left(|v|^{2}v-|w|^{2}w \right)e^{2 Y} \overline{r} \e^{2 Y_N}\d x\right|.
    \end{align*}
    Using \cref{product_rule,product-rule-W(-k;r)-Wkp,action_D&x/Wsp,lem:product-with-exp(aY)}, it gives
    \begin{align*}
        \left|\frac{1}{2}\frac{\d}{\d t} \int_{\mathbb{R}^2} |r|^2 \e^{2Y_N}\d x\right|
        &\lesssim\left(\left|\nabla(Y-Y_N)\right|_{\W^{-\kappa,\infty}_x}+\left|x(Y-Y_N)\right|_{\W^{-\kappa,\infty}_x}\right)\left|r\right|_{\W^{1+\kappa,2}_x}\left|r\right|_{\W^{\kappa,2}_x}\left|\e^{2Y_N}\right|_{1+\W^{\kappa,\infty}_x}\\ 
        &+\left|:|\nabla Y|^2:-:|\nabla Y_N|^2:\right|_{\W^{-\kappa,\infty}_x}\left|r\right|_{\W^{\kappa,2}_x}^2\left|\e^{2Y_N}\right|_{1+\W^{\kappa,\infty}_x}\\ 
        &+ \left(|v|^{2}_{\mathbb{L}^\infty_{x}}+|w|^{2}_{\mathbb{L}^\infty_{x}}\right)\left|e^{2 Y}\right|_{\mathbb{L}^\infty_{x}} \int_{\mathbb{R}^2} |r|^2\e^{2 Y_N}\d x\\ 
        &\lesssim\left(\left|\nabla(Y-Y_N)\right|_{\W^{-\kappa,\infty}_x}+\left|x(Y-Y_N)\right|_{\W^{-\kappa,\infty}_x}\right)\left|r\right|_{\mathbb{L}^\infty_{t}\W^{\sigma,2}_x}^2\left|\e^{2Y_N}\right|_{1+\W^{\kappa,\infty}_x}\\ 
        &+\left|:|\nabla Y|^2:-:|\nabla Y_N|^2:\right|_{\W^{-\kappa,\infty}_x}\left|r\right|_{\mathbb{L}^\infty_{t}\W^{\sigma,2}_x}^2\left|\e^{2Y_N}\right|_{1+\W^{\kappa,\infty}_x}\\ 
        &+ \left(|v|^{2}_{\mathbb{L}^\infty_{t}\W^{\sigma,2}_x}+|w|^{2}_{\mathbb{L}^\infty_{t}\W^{\sigma,2}_x}\right)\left|e^{2 Y}\right|_{\mathbb{L}^\infty_{x}} \int_{\mathbb{R}^2} |r|^2\e^{2 Y_N}\d x.
    \end{align*}
    Recall that $\delta \in\left(0,\frac{2\kappa}{3}\right)$, so there exists a positive random variable $C$ finite on $\Omega_0$ (by definition of $\Omega_0$) such that for any $N\in\N$, it holds
    $$\frac{\d}{\d t} \int_{\mathbb{R}^2} |r|^2\e^{2 Y_N}\d x \leqslant C\left(\left|r\right|_{\mathbb{L}^\infty_{t}\W^{\sigma,2}_x}^2 \lambda_N^{-\delta}+ \left(|v|^{2}_{\mathbb{L}^\infty_{t}\W^{\sigma,2}_x}+|w|^{2}_{\mathbb{L}^\infty_{t}\W^{\sigma,2}_x}\right)\left|e^{2 Y}\right|_{\mathbb{L}^\infty_{x}} \int_{\mathbb{R}^2} |r|^2\e^{2 Y_N}\right).$$
    Then by Gronwall's lemma, we obtain
    $$\forall N\in\N,\;\sup_{t\in[a,b]}\int_{\mathbb{R}^2} |r|^2\e^{2 Y_N}\d x \leqslant C\left|r\right|_{\mathbb{L}^\infty_{t}\W^{\sigma,2}_x}^2 \lambda_N^{-\delta} \e^{(b-a) C \left(|v|^{2}_{\mathbb{L}^\infty_{t}\W^{\sigma,2}_x}+|w|^{2}_{\mathbb{L}^\infty_{t}\W^{\sigma,2}_x}\right)\left|e^{2 Y}\right|_{\mathbb{L}^\infty_{x}}}$$
    which in addition to \cref{Equation_unicity_control} implies that $v=w$ almost surely.\endproof
\end{demo}

We will only prove the existence of a solution in the defocusing case (i.e. \cref{Theorem-cubic-defocusing}). In view of \cref{Proposition_energy_bound_focusing}, analogous proofs will hold on the conditional probability space $(\overline{\Omega}_L,\mathcal{A}_{|\overline{\Omega}_L},\P_{|\overline{\Omega}_L})$ in the focusing case if we restrict ourselves to initial data $v_0\in\W^{2,2}$ satisfying $|v_0|_{\mathbb{L}^2_{x}}\leqslant L$. So, assume $\lambda\leqslant0$. For $N\in\N$, denote by $v^N$ the unique global solution to \cref{Eq:2D-transformed-renorm-regularized}. First as in the linear case, we obtain a diverging bound but this time on $|v^N|_{\mathbb{L}^\infty_t\W^{\sigma,2}_x}$ for any $\sigma\in\left(\frac{3}{2},2\right)$. In order to obtain this bound, we prove the following generalized Brezis-Gallouet inequality (for the classical one see \cite{BREZIS1980677}).

\begin{prop}{Generalized Brezis-Gallouet inequality}\label{Generalized_Brezis-Gallouet_inequality}

    Let $\sigma>1$, there exists a constant $C>0$ such that 
    $$\forall v\in\W^{\sigma,2},\; |v|_{\mathbb{L}^\infty_{x}}\leqslant C \left(1+ |v|_\Sigma\sqrt{1+ \ln(1+|v|_{\W^{\sigma,2}_x})}\right)$$
    
\end{prop}

\begin{demo}
    By continuity of both sides and density of $\mathcal{S}(\mathbb{R}^2)$ in $\W^{\sigma,2}$, it is sufficient to obtain the claim for $v\in\mathcal{S}(\mathbb{R}^2)$. Then, we have for any $R>0$:
    \begin{align*}
        |v|_{\mathbb{L}^\infty_{\mathbb{R}^2}} &\lesssim \int_{|\eta|<R} \langle \eta\rangle \mathcal{F} v(\eta) \frac{\d \eta}{\langle \eta\rangle} + \int_{|\eta|\geqslant R} \langle \eta\rangle^\sigma \mathcal{F} v(\eta) \frac{\d \eta}{\langle \eta\rangle^\sigma}\\ 
        &\lesssim |v|_{\Sigma}\sqrt{\ln(1+R)} + \frac{|v|_{\W^{\sigma,2}}}{(1+R^2)^{\frac{\sigma-1}{2}}}\\ 
        &\lesssim |v|_{\Sigma}\sqrt{\ln(1+R)} + \frac{|v|_{\W^{\sigma,2}}}{1+R^{\sigma-1}}
    \end{align*}
    using Cauchy-Schwartz inequality. We conclude choosing $R=\left(\frac{|v|_{\W^{\sigma,2}}}{1+|v|_\Sigma}\right)^{\frac{1}{\sigma-1}}\leqslant |v|_{\W^{\sigma,2}}^{\frac{1}{\sigma-1}}$ and using properties of the logarithm.\endproof
\end{demo}

Using this inequality, we obtain the following diverging bound which grows polynomially in $N$.

\begin{prop}\label{Prop:diverging-bound-W^(sigma;2)-defocusing}
    Almost surely, for all $v_0\in\W^{2,2}$, for all $a<0<b$, for all $\sigma\in(1,2)$ and all $\kappa\in(0,1)$, there exists $C>0$ such that
    $$\forall N\in\N,\; \left|v^N\right|_{\mathbb{L}^\infty_{[a,b]}\W^{\sigma,2}}\leqslant C \lambda_N^{\kappa}.$$
\end{prop}

Then, as in the linear case, we obtain a Cauchy estimate in $\mathbb{L}^2(\mathbb{R}^2)$.

\begin{prop}\label{Proposition-Bound-L2-Diff-cubic-defocusing}

    Almost surely, for all $v_0\in\W^{2,2}$, for all $a<0<b$, for all $\kappa\in(0,1)$ and all $\delta\in\left(0,\frac{1}{3}\right)$, there exists $C>0$ such that
    $$\forall N>M\geqslant 0,\; \left|v^N-v^M\right|_{\mathbb{L}^\infty_{[a,b]}\mathbb{L}^2_{\mathbb{R}^2}} \leqslant C \lambda_N^{\kappa}\lambda_M^{-\delta}.$$
    
\end{prop}

The proofs of both \cref{Prop:diverging-bound-W^(sigma;2)-defocusing,Proposition-Bound-L2-Diff-cubic-defocusing} are given in \cref{appendix_proof_bounds_cubic}. We can now prove the theorem.

\begin{demo}{of \cref{Theorem-cubic-defocusing}.}

    Let $\overline{\Omega}=\Omega'\cap\Omega''\cap\Omega_0$ where $\Omega'$, $\Omega''$ and $\Omega_0$ are the full probability event of \cref{Prop:diverging-bound-W^(sigma;2)-defocusing,,Proposition-Bound-L2-Diff-cubic-defocusing,Proposition-unicity-nonlinear-equation} respectively. Then $\P(\overline{\Omega})=1$ and for any $N\in\N$ and any $v_0\in\W^{2,2}$ and any $\omega\in\overline{\Omega}$, $v^N(\omega)$ is global.\\ 

    Let $\omega\in\overline{\Omega}$, $v_0\in\W^{2,2}$, $a<0<b$ and $\sigma,\sigma'\in(1,2)$ such that $\sigma'>\sigma$, $\delta\in\left(0,\frac{1}{3}\right)$ and let $\kappa = \left(1-\frac{\sigma}{\sigma'}\right)\frac{\delta}{2}\in(0,1)$. Then by interpolation, for any $N>M\geqslant 0$, it holds
    $$\left|v^N(\omega)-v^M(\omega)\right|_{\mathbb{L}^\infty_{t}\W^{\sigma,2}_x} \leqslant \left|v^N(\omega)-v^M(\omega)\right|_{\mathbb{L}^\infty_{t}\W^{\sigma',2}_x}^{\frac{\sigma}{\sigma'}}\left|v^N(\omega)-v^M(\omega)\right|_{\mathbb{L}^\infty_{t}\mathbb{L}^2_{x}}^{1-\frac{\sigma}{\sigma'}}.$$
    By \cref{Prop:diverging-bound-W^(sigma;2)-defocusing,Proposition-Bound-L2-Diff-cubic-defocusing}, there exists $C(\omega,v_0)\in(0,+\infty)$  (which is independent of $N$ and $M$, and may vary from line to line) such that
    \begin{equation}
        \forall N>M\geqslant 0,\; \left|v^N(\omega)-v^M(\omega)\right|_{\mathbb{L}^\infty_{t}\W^{\sigma,2}_x} \leqslant C(\omega,v_0)\lambda_N^\kappa\lambda_M^{-(1-\frac{\sigma}{\sigma'})\delta}.\label{Equation_proof_cubic_theorem_interpolated}
    \end{equation}
    Let $n\in\N$, $M=2^n$ and $N=2^{n+1}$. Then \cref{Equation_proof_cubic_theorem_interpolated} gives
    $$\left|v^{2^{n+1}}(\omega)-v^{2^n}(\omega)\right|_{\mathbb{L}^\infty_{t}\W^{\sigma,2}_x} \leqslant C(\omega,v_0)2^{\frac{n}{2}\left(\kappa-(1-\frac{\sigma}{\sigma'})\delta\right)}=C(\omega,v_0) 2^{-n(1-\frac{\sigma}{\sigma'})\frac{\delta}{4}},$$
    which is summable. Hence the sequence $\left(v^{2^n}(\omega)\right)_{n\in\N}$ is Cauchy in the Banach space $\mathcal{C}([a,b],\W^{\sigma,2})$ and thus converges to an element $v(\omega)$. As in the linear case, using
    $$\left|v^{N}(\omega)-v(\omega)\right|_{\mathbb{L}^\infty_{t}\W^{\sigma,2}_x}\leqslant \left|v^{N}(\omega)-v^{2^{\lfloor \log_2 N\rfloor}}(\omega)\right|_{\mathbb{L}^\infty_{t}\W^{\sigma,2}_x}+\left|v^{2^{\lfloor \log_2 N\rfloor}}(\omega)-v(\omega)\right|_{\mathbb{L}^\infty_{t}\W^{\sigma,2}_x}$$
    in addition to \cref{Equation_proof_cubic_theorem_interpolated} with a well-chosen $\kappa$ shows the convergence of $v^N(\omega)$ to $v(\omega)$ in $\mathcal{C}([a,b],\W^{\sigma,2})$. The function $v$ is measurable as an almost sure limit of measurable functions. Using \cref{product-rule-W(-k;r)-Wkp,Cor:convergence-Y_N-W(1-kappa;infty),Lem:Conv-W(1-kappa;inf)-exp(aYN),Cor:cauchy-criter-W(-kappa;infty)}, it is not difficult to show that $-\mathrm{i}\p_t v^{N}(\omega)$ converges in $\mathcal{C}([a,b],\W^{\sigma-2,2})$ to $$Hv(\omega)+2\nabla v(\omega)\cdot \nabla Y(\omega) + xv(\omega)\cdot xY(\omega)+\wick{|\nabla Y|^2}(\omega) v(\omega)+\lambda |v(\omega)|^{2} v(\omega) \e^{2 Y(\omega)} $$ so that $v\in\mathcal{C}^1([a,b],\W^{\sigma-2,2})\cap\mathcal{C}([a,b],\W^{\sigma,2})$ a.s. and is solution of \cref{Eq:2D-transformed-renorm}. Moreover, as in the linear case, the regularity of $v$ is sufficient to pass to the limit conservation laws.\endproof

\end{demo}

\begin{rem}{~}

    \begin{itemize}
        \item In the case of subcubic nonlinearities, $|u|^{2\gamma}u$ with $\gamma\in\left(\frac{1}{2},1\right)$, analogous of \cref{Proposition_energy_bound_defocusing,Proposition-unicity-nonlinear-equation,Prop:diverging-bound-W^(sigma;2)-defocusing,Proposition-Bound-L2-Diff-cubic-defocusing,Theorem-cubic-defocusing} holds and there exists almost surely, for all initial data in $\W^{2,2}$, a unique solution which lives in $\mathcal{C}(\R,\W^{2-,2})$.
        \item If instead of working with deterministic initial data in the focusing case, one choose a random initial datum $v_0\in\W^{2,2}$ a.s. such that 
        $$\lambda |v_0|_{\mathbb{L}^2_{x}}^2 \left|\e^{-2Y}\right|_{\mathbb{L}^\infty_{x}}^2\left|\e^{2Y}\right|_{\mathbb{L}^\infty_{x}}\left|\e^{4Y}\right|_{\mathbb{L}^\infty_{x}}<4 $$
        holds almost surely, then by a slight modification of previous arguments, one can show that there exists almost a unique global solution.
    \end{itemize}
    
\end{rem}

\appendix

\section{Proofs of results from Section \ref{section_preliminaries}}\label{appendix_proof_preli}

\subsection{Hermite-Sobolev spaces}

\begin{demo}{of \cref{Proposition_creation_annhil_Wsp}}

    It is known from \cite{BongioanniHSspaces} that, for $k\in\N$, an equivalent norm on $\W^{k,p}$ is given by 
    $$|u|_{k,p} = |u|_{\mathbb{L}^p_x}+\sum_{m=1}^k \sum_{1\leqslant|i_1|,\dotsc,|i_m|\leqslant 2} |A_{i_1}\cdots A_{i_m} u|_{\mathbb{L}^p_x}.$$
    Then, for $u\in\W^{k,p}$ with $k\in\N^*$ and $i\in\{\pm 1,\pm 2\}$, we have 
    $$|A_i u|_{k-1,p} = |A_iu|_{\mathbb{L}^p_x}+\sum_{m=1}^{k-1} \sum_{1\leqslant|i_1|,\dotsc,|i_m|\leqslant 2} |A_{i_1}\cdots A_{i_m} A_i u|_{\mathbb{L}^p_x} \leqslant |u|_{k,p},$$
    so $A_i\in\mathcal{L}\left(\W^{k,p},\W^{k-1,p}\right)$. The case $k\in-\N$ is shown by duality and the general case $s\in\R$ is obtained by interpolation.\endproof
     
\end{demo}

\begin{demo}{of \cref{Sobolev-embeddings}.}

    The first point follows from usual Sobolev embedding of $W^{s-\sigma,p}$ in $\mathbb{L}^q(\mathbb{R}^2)$ and continuous embedding of $\W^{s-\sigma,p}$ in $W^{s-\sigma,p}$. When the condition is sharp, one can use the compact embedding of $\W^{\eps,p}$ in $\mathbb{L}^p(\mathbb{R}^2)$ (for $\eps>0$ and $p\in(1,+\infty)$) to obtain compactness of the embedding.\\ 

    It remains to prove continuous embedding in the case $p<q=+\infty$. Assume first that $\sigma<0$ and let $p'$ be such that $\frac{1}{p}+\frac{1}{p'}=1$. As $\frac{1}{p}-\frac{s}{2}<\frac{1}{q}-\frac{\sigma}{2}$, it follows $1<\frac{1}{p'}+\frac{s-\sigma}{2}$. For any $k\in\N$, $\W^{k,1}$ is continuously embedded in $W^{k,1}$ which is continuously embedded in $W^{k+\sigma-s,p'}$ by classical Sobolev embeddings. By interpolation, $\W^{-\sigma,1}$ is continuously embedded in $W^{-s,p'}$. By duality and continuous embedding of $\W^{s,p}$ in $W^{s,p}$ we obtain the claim for $\sigma<0$.\\ 

     Now, for $\sigma=2k\in2\N$, the claim follows from usual continuous embeddings of $W^{s,p}$ in $W^{2k,\infty}$ for the Sobolev part and using the continuous injection of $\W^{s-2k,p}$ in $\mathbb{L}^\infty(\mathbb{R}^2)$ combined with \cref{Corollary_action_D_and_x_on_Wsp} to obtain that $\langle x\rangle^{2k}:\W^{s,p} \to \mathbb{L}^\infty(\mathbb{R}^2)$ is continuous. The general case $\sigma\geqslant0$ follows by interpolation.\endproof
\end{demo}

\subsection{The smooth truncation}

\begin{demo}{of \cref{lemme-low-freq-estim}.}

    As $S_N$ commute with $(-H)^{\frac{\alpha}{2}}$, it is sufficient to prove the case $\alpha=0$. By interpolation, it is sufficient to only prove the case $s=2n$ for $n\in\N$. Let $N\in\N$ and $\phi\in\mathbb{L}^p(\mathbb{R}^2)$, write it $\phi=\sum_{k\in\N^2}\phi_k h_k$, then we have 
    $$|S_N\phi|_{\W^{2n,p}_x}^p = \int_{\mathbb{R}^2}\left|\sum_{k\in\N^2} \chi_{k,N} \lambda_k^{2n} \phi_k\right|^p\d x= \lambda_N^{2np}\int_{\mathbb{R}^2}\left|\sum_{k\in\N^2} \chi_{k,N} \left(\frac{\lambda_k^2}{\lambda_N^2}\right)^n \phi_k h_k(x)\right|^p\d x.$$
    Define $\psi_n(\lambda)=\lambda^n\chi(\lambda)$, then $\psi\in\mathcal{S}(\R)$, so
    by \cref{Lemma_psi(-thetaH)}, it follows
    $$|S_N\phi|_{\W^{2n,p}_x}^p = \lambda_N^{2np}\left|\psi_n\left(\frac{-H}{\lambda_N^2}\right)\phi\right|_{\mathbb{L}^p_{x}}^p\lesssim \lambda_N^{2np}\left|\phi\right|_{\mathbb{L}^p_{x}}^p.$$
    \endproof
\end{demo}

\begin{demo}{of \cref{lemme-high-freq-estim-Ws2}.}

    Let $\phi=\sum_{k\in\N^2}\phi_k h_k\in\W^{\alpha+s,2}$. It holds
    $$|\phi-S_N\phi|_{\W^{\alpha,2}_x}^2 = \sum_{k\in\N^2} (1-\chi_{k,N}^2) \lambda_k^{2\alpha} |\phi_k|^2 \leqslant \lambda_{\left\lfloor\frac{N}{2}\right\rfloor}^{-2s} \sum_{k\in\N^2} \lambda_k^{2\alpha+2s} |\phi_k|^2 \leqslant \lambda_{\left\lfloor\frac{N}{2}\right\rfloor}^{-2s}|\phi|_{\W^{\alpha+s,2}_x}^2.$$

    \endproof
\end{demo}

\section{Proofs of results from Section \ref{Section_noise}}\label{appendix_proof_proba}

\subsection{Noise regularity}

\begin{demo}{of \cref{Lemma_bound_from_moment}.}

    Let $\beta<\alpha$ and $p>\max\left(1,\frac{2}{\alpha-\beta}\right)$. By Markov's inequality, we have
    $$\P\left(\lambda_N^\beta |X_N|>1\right) \leqslant \lambda_N^{-p(\alpha-\beta)}\left(\sup_{N\in\N} \lambda_N^\alpha |X_N|_{\mathbb{L}^p_\omega}\right)^p.$$

    Let $A_N = \left\{\lambda_N^\beta |X_N|>1\right\}$. The previous inequality shows $\sum_{N\in\N} \P(A_N)<+\infty$, thus, by the Borel-Cantelli lemma $\displaystyle\P\left(\limsup_{N\to+\infty} A_N\right)=0.$
    Hence, almost surely, there exists $n\in\N$ such that for any $N\geqslant n$, $\lambda_N^\beta |X_N|\leqslant 1$. This implies that almost surely $\displaystyle\sup_{N\in\N} \lambda_N^\beta |X_N| <+\infty.$\endproof
    
\end{demo}

\begin{demo}{of \cref{Lem:convergence-Y_N-W(1-s;q)}.}

    By Hölder's inequality, it is sufficient to prove the case $p\geqslant q$. By Minkowski's inequality and gaussianity, we have 
    $$\E\left[|Y-Y_N|^p_{\W^{1-s,q}_x}\right]^{\frac{q}{p}} \lesssim \int_{\mathbb{R}^2} \E\left[\left|\sum_{k\in\N^2} (1-\chi_{k,N}) \lambda_k^{-(1+s)} h_k(x) \xi_k\right|^2\right]^{\frac{q}{2}} \d x\lesssim \lambda_{N}^{-q(s-s')} |K|^q_{\mathbb{L}^q_x\W^{1-s',2}_y}.$$
    The assumption $s'q>2$ implies $1-s'<1-\frac{2}{q}$ and we conclude using \cref{regu_kernel_-H-1}.\endproof
    
\end{demo}

\begin{demo}{of \cref{Cor:convergence-Y_N-W(1-kappa;infty)}.}

    Let $\epsilon = \sqrt{\frac{\kappa'}{\kappa}}$, $s=\epsilon \kappa$, $s'=(1-\epsilon)s$ and $q>\max\left(\frac{2}{s'},\frac{2}{\kappa-s}\right)$. By Sobolev embeddings and \cref{action_D&x/Wsp}, there exists $c>0$ such that
    $$\max\left(|Y-Y_N|_{\W^{1-\kappa,\infty}_x},|\nabla Y-\nabla Y_N|_{\W^{-\kappa,\infty}_x},|xY-xY_N|_{\W^{-\kappa,\infty}_x} \right)\leqslant c |Y-Y_N|_{\W^{1-s,q}_x}.$$
    Denote by $C_N = c |Y-Y_N|_{\W^{1-s,q}_x} \lambda_N^{s-s'}$. As $s-s'=\kappa'$, we obtain
    $$\max\left(|Y-Y_N|_{\W^{1-\kappa,\infty}_x},|\nabla Y-\nabla Y_N|_{\W^{-\kappa,\infty}_x},|xY-xY_N|_{\W^{-\kappa,\infty}_x} \right)\leqslant C_N \lambda_N^{-\kappa'}.$$
    Let $p\in[1,+\infty)$, by \cref{Lem:convergence-Y_N-W(1-s;q)},  $\sup_{N\in\N} |C_N|_{\mathbb{L}^p_\omega}<+\infty.$
    The only point remaining is to show that almost surely $\sup_{N\in\N} C_N<+\infty$. Let $\kappa''\in(\kappa',\kappa)$, then by taking $s=\sqrt{\kappa\kappa'}$, $s'=\left(1-\sqrt{\frac{\kappa''}{\kappa}}\right)s$ and $q>\max\left(\frac{2}{s'},\frac{2}{\kappa-s}\right)$, we can conclude using \cref{Lem:convergence-Y_N-W(1-s;q),Lemma_bound_from_moment}.\endproof
    
\end{demo}

\begin{demo}{of \cref{Corollary_divergence_moment_norm_Linfty_GradYN_xYN}.}

    Let $q>\frac{4}{\kappa}$ and $s=\frac{\kappa}{2}$, by Sobolev embedding and \cref{Corollary_action_D_and_x_on_Wsp}, there exists a constant $c>0$ such that $$\max\left(|\nabla Y_N|_{\mathbb{L}^{\infty}_{x}},|x Y_N|_{\mathbb{L}^{\infty}_{x}}\right)\leqslant c |Y_N|_{\W^{1+s,q}_x}.$$
    Hence, using \cref{lemme-low-freq-estim,bound-S_N-Wsq}, there exists $c'>0$ such that
    $$\max\left(|\nabla Y_N|_{\mathbb{L}^{\infty}_{x}},|x Y_N|_{\mathbb{L}^{\infty}_{x}}\right)\leqslant c' |Y|_{\W^{1-s,q}_x} \lambda_N^{\kappa},$$
    where we used that $\kappa=2s$.\endproof
    
\end{demo}

\begin{demo}{of \cref{Lem:moment-boud-expYN-Wsq}.}

    Due to the condition on $s$, $\W^{1-s,q}$ is continuously embedded in $\mathbb{L}^\infty(\mathbb{R}^2)$. Thus, $\W^{1-s,q}$ is an algebra and there exists a constant $c\geqslant1$ such that for all $u,v\in\W^{s,q}$ we have
    \begin{equation}
        |uv|_{\W^{1-s,q}_x}\leqslant c |u|_{\W^{1-s,q}_x}|v|_{\W^{1-s,q}_x}.\label{Eq:algebra-property}
    \end{equation}
    Using Minkowski's inequality, we obtain
    $$\forall N\in\N,\;\E\left[\left|\e^{a Y_N}-1\right|_{\W^{1-s,q}_x}^p\right]\leqslant\E\left[\e^{cp|a| \left|Y_N\right|_{\W^{1-s,q}_x}}\right]\text{ and }\E\left[\left|\e^{a Y}-1\right|_{\W^{1-s,q}_x}^p\right]\leqslant\E\left[\e^{c p|a| \left|Y\right|_{\W^{1-s,q}_x}}\right].$$
    By \cref{bound-S_N-Wsq}, there exists $c'\geqslant 1$ such that, for any $N\in\N$, we have $\left|Y_N\right|_{\W^{1-s,q}_x}\leqslant c'|Y|_{\W^{1-s,q}_x}.$
    By Fernique's theorem (see \cite{MR266263}), there exists $\beta>0$ such that $\E\left[\e^{\beta |Y|_{\W^{1-s,q}_x}^2}\right]<+\infty.$
    So, let $c''=\frac{(c c'p|a|)^2}{\beta}$, then Young's inequality implies,
    $$\forall N\in\N,\;\E\left[\left|\e^{a Y_N}-1\right|_{\W^{1-s,q}_x}^p\right]\leqslant\e^{c''}\E\left[\e^{\beta |Y|_{\W^{1-s,q}_x}^2}\right]$$ and $$\E\left[\left|\e^{a Y}-1\right|_{\W^{1-s,q}_x}^p\right]\leqslant\e^{c''}\E\left[\e^{\beta |Y|_{\W^{1-s,q}_x}^2}\right].$$\endproof

\end{demo}

\begin{demo}{of \cref{Lem:Conv-W(1-kappa;inf)-exp(aYN)}.}

    The second claim follows easily from the first one using \cref{Lemma_bound_from_moment}. In order to prove the first claim, let $q>\frac{4}{\kappa-\kappa'}$ and $s\in\left(\kappa'+\frac{2}{q},\kappa-\frac{2}{q}\right)$. Then $\W^{1-s,q}$ is continuously embedded in $\W^{1-\kappa,\infty}$ and is an algebra, thus \cref{Eq:algebra-property} is verified. Denote by $C$ the constant of Sobolev embedding of $\W^{1-s,q}$ in $\W^{1-\kappa,\infty}$. Using the algebra property of $\W^{1-s,q}$ norms, we obtain
    $$\left|\e^{aY_N}-\e^{aY}\right|_{\W^{1-\kappa,\infty}_x}\leqslant cC\left|\e^{aY}\right|_{\W^{1-s,q}_x}\left|\e^{a(Y_N-Y)}-1\right|_{\W^{1-s,q}_x} .$$
    By Minkowski's inequality, we have
    \begin{equation*}
        \left|\e^{a(Y_N-Y)}-1\right|_{\W^{1-s,q}_x} \leqslant \sum_{n=1}^{+\infty} \frac{c^{n-1}|a|^n|Y-Y_N|_{\W^{1-s,q}_x}^n}{n!} \leqslant |a||Y-Y_N|_{\W^{1-s,q}_x}\e^{c|a||Y-Y_N|_{\W^{1-s,q}_x}}.
    \end{equation*}
    Using Hölder's inequality, \cref{bound-S_N-Wsq,Lem:moment-boud-expYN-Wsq,Lem:convergence-Y_N-W(1-s;q)} with $s'=s-\kappa'$ we obtain the first claim. \endproof
\end{demo}

\subsection{Renormalization}

\begin{demo}{of \cref{nelson-estimate-norm-W-alpha-q}.}

    By Minkowski's inequality and \cref{nelson-estimate}, we have
    $$\left|:|\nabla Y_N|^2:\right|_{\mathbb{L}^{q'}_\omega\W^{\alpha,q}_x}^q\leqslant \int_{\mathbb{R}^2}  \E\left[\left|(-H)^{\frac{\alpha}{2}}:|\nabla Y_N|^2:(x)\right|^{q'}\right]^{\frac{q}{q'}}\d x\leqslant (q'-1) \left|(-H)^{\frac{\alpha}{2}}:|\nabla Y_N|^2:\right|_{\mathbb{L}^{q}_x\mathbb{L}^2_\omega}.$$
    Using one more time Minkowski's and Hölder's inequalities, we obtain the claim.\endproof
    
\end{demo}

\begin{demo}{of \cref{Cor:cauchy-criter-W(-kappa;infty)}.}

    First, we show the moment estimates in the case $q<+\infty$. Recall that \cref{conv_renorm} implies
    $$\forall N\geqslant M,\; \left|\wick{|\nabla Y_N|^2}-\wick{|\nabla Y_M|^2}\right|_{\mathbb{L}^q_\omega\W^{-s,q}_x} \lesssim \lambda_M^{-\delta}.$$
    Nelson's estimate and Hölder's inequality then implies
    $$\forall N\geqslant M,\; \left|\wick{|\nabla Y_N|^2}-\wick{|\nabla Y_M|^2}\right|_{\mathbb{L}^p_\omega\W^{-s,q}_x} \lesssim \lambda_M^{-\delta}.$$

    In the case $q=+\infty$, let $q'>\max\left(\frac{5}{s},\frac{10}{2s-3\delta}\right)$ and define $s'=s-\frac{3}{q'}$. By definition of $q'$, we have $s'>\frac{2}{q'}$ and the Sobolev embedding of $\W^{-s',q'}$ in $\W^{-s,q}$. Moreover, $\delta<\frac{2s'}{3}-\frac{4}{3q'}$. Hence the previous estimate leads to
    $$\forall N\geqslant M,\; \left|\wick{|\nabla Y_N|^2}-\wick{|\nabla Y_M|^2}\right|_{\mathbb{L}^p_\omega\W^{-s,q}_x} \lesssim \left|\wick{|\nabla Y_N|^2}-\wick{|\nabla Y_M|^2}\right|_{\mathbb{L}^p_\omega\W^{-s',q'}_x} \lesssim \lambda_M^{-\delta}.$$
    Let $M\in\N$, by triangular inequality, we have
    $$\left|\wick{|\nabla Y_M|^2}\right|_{\mathbb{L}^p_\omega\W^{-s,q}_x} \leqslant \left|\wick{|\nabla Y_M|^2}-\wick{|\nabla Y|^2}\right|_{\mathbb{L}^p_\omega\W^{-s,q}_x} + \left|\wick{|\nabla Y|^2}\right|_{\mathbb{L}^p_\omega\W^{-s,q}_x}.$$
     \Cref{nelson-estimate-norm-W-alpha-q} shows that $\left|\wick{|\nabla Y|^2}\right|_{\mathbb{L}^p_\omega\W^{-s,q}_x}<+\infty$ and the previous point shows that (by sending $N$ to $+\infty$) 
     \begin{equation}
         \forall M\in\N,\;\left|\wick{|\nabla Y_M|^2}-\wick{|\nabla Y|^2}\right|_{\mathbb{L}^p_\omega\W^{-s,q}_x}\lesssim \lambda_M^{-\delta} \lesssim 1.\label{Proof-cauchy-criter-Wick-equation}
     \end{equation}
    Hence $$\max\left(\left|\wick{|\nabla Y|^2}\right|_{\mathbb{L}^p_\omega\W^{-s,q}_x},\sup_{N\in\N}\left|\wick{|\nabla Y_N|^2}\right|_{\mathbb{L}^p_\omega\W^{-s,q}_x}\right)<+\infty$$ and we conclude using \cref{Lemma_bound_from_moment,Proof-cauchy-criter-Wick-equation}.\endproof
    
\end{demo}

\begin{demo}{of \cref{control-norm-Lq-Wick-product}.}

    By Fubini's theorem and Nelson's estimate, we have 
    $$\E\left[\left|\wick{|\nabla Y_N|^2}\right|_{\mathbb{L}^{q}_{x}}^q\right] = \int_{\mathbb{R}^2}  \E\left[\left|\wick{|\nabla Y_N|^2}(x)\right|^q\right] \d x\lesssim \int_{\mathbb{R}^2}  \E\left[\left|\wick{|\nabla Y_N|^2}(x)\right|^2\right]^{\frac{q}{2}} \d x.$$
    Using $\wick{|\nabla Y_N|^2}(x) = |\nabla Y_N(x)|^2-\E\left[|\nabla Y_N(x)|^2\right]$, one has 
    $$\E\left[\left|\wick{|\nabla Y_N|^2}(x)\right|^2\right] = \var\left(|\nabla Y_N(x)|^2\right)\leqslant \E\left[\left|\nabla Y_N(x)\right|^4\right].$$
    Using Jensen's inequality, we obtain
    $$\E\left[\left|\wick{|\nabla Y_N|^2}\right|_{\mathbb{L}^q_{x}}^q\right]\lesssim \int_{\mathbb{R}^2}  \E\left[\left||\nabla Y_N|(x)\right|^4\right]^{\frac{q}{2}} \d x\lesssim \int_{\mathbb{R}^2}  \E\left[\left|\left(|\nabla Y_N|^2\right)(x)\right|^q\right] \d x\lesssim \left||\nabla Y_N|^2\right|_{\mathbb{L}^{q}_\omega\mathbb{L}^{q}_{x}}^{q},$$
    and we conclude using Hölder's inequality. \endproof
    
\end{demo}

\begin{demo}{of \cref{estim-Lq-norm-Wick-product}.}

    Using Hölder's inequality for $p< q$ and \cref{nelson-estimate-norm-W-alpha-q} for $p>q$, it is sufficient to prove the case $p=q$. Using \cref{control-norm-Lq-Wick-product,lemme-low-freq-estim}, we have 
    $$ \left|\wick{|\nabla Y_N|^2}\right|_{\mathbb{L}^q_\omega\mathbb{L}^q_{x}}\lesssim\left|\nabla Y_N\right|_{\mathbb{L}^{2q}_\omega\mathbb{L}^{2q}_{x}}^2\lesssim \lambda_N^{s} |\nabla Y_N|_{\mathbb{L}^{2q}_\omega\W^{-\frac{s}{2},2q}_x}^{2}$$
    and we conclude using \cref{Lem:convergence-Y_N-W(1-s;q)}.\endproof
    
\end{demo}

\subsection{The regularized equation}

\begin{demo}{of \cref{Proposition_energy_bound_linear}}

First using definition of the transformed energy, we obtain
\begin{align*}
    \frac{1}{2}|v^N(t)|^2_\Sigma &= \frac{1}{2}\int_{\mathbb{R}^2} \left(|\nabla v^N(t)|^2+|xv^N(t)|^2\right)\e^{2Y_N}\e^{-2Y_N} \d x\\ 
    &\leqslant |\e^{-2Y_N}|_{\mathbb{L}^\infty_x}\left(\Tilde{\mathcal{E}}_N(v^N(t))+\int_{\mathbb{R}^2} |xv^N(t)|^2Y_N\e^{2Y_N}\d x+\int_{\mathbb{R}^2} |v^N(t)|^2:|\nabla Y_N|^2:\e^{2Y_N}\d x\right).
\end{align*}

Conservation of the transformed energy implies 
\begin{equation}
    \begin{split}
        \frac{1}{2}|v^N(t)|^2_\Sigma&\leqslant |\e^{-2Y_N}|_{\mathbb{L}^\infty_x}\Tilde{\mathcal{E}}_N(v_0)+|\e^{-2Y_N}|_{\mathbb{L}^\infty_x}|\e^{2Y_N}|_{\mathbb{L}^\infty_x}\int_{\mathbb{R}^2} |xv^N(t)|^2|Y_N|\d x\\ 
    &+\left||\e^{-2Y_N}|_{\mathbb{L}^\infty_x}\int_{\mathbb{R}^2} |v^N(t)|^2:|\nabla Y_N|^2:\e^{2Y_N}\d x\right|.
    \end{split}\label{Equation_bound_Sigma_from_E}
\end{equation}

We treat the two integral terms separately.

\begin{itemize}
    \item Let $q>2$, $s>0$ such that $qs>2$ and $q'>1$ verifying $\frac{1}{q}+\frac{1}{q'}=1$. Using Hölder's inequality, Sobolev embeddings and interpolation, we obtain
        \begin{align*}
            \int_{\mathbb{R}^2} |xv^N(t)|^2 |Y_N|\d x &\lesssim \left||x|^{\frac{1+s}{2}}v^N(t)\right|_{\mathbb{L}^{2q'}_x}^2\left||x|^{1-s}Y_N\right|_{\mathbb{L}^{q}_x}\\ 
            &\lesssim \left|v^N(t)\right|_{\W^{\frac{1+s}{2},2q'}_x}^2\left|Y_N\right|_{\W^{1-s,q}_x}\\ 
            &\lesssim \left|v^N(t)\right|_{\mathbb{L}^2_x}^{1-\frac{2}{q}-s}\left|v^N(t)\right|_{\Sigma}^{1+\frac{2}{q}+s}\left|Y_N\right|_{\W^{1-s,q}_x}.
        \end{align*}

        Let $q_1>1$ such that $q_1\left(1+\frac{2}{q}+s\right)=2$ and $p_1$ such that $\frac{1}{p_1}+\frac{1}{q_1}=1$. Using Young's inequality and mass conservation, there exists a constant $c_1>0$ such that
        \begin{equation*}
            |\e^{-2Y_N}|_{\mathbb{L}^\infty_x}|\e^{2Y_N}|_{\mathbb{L}^\infty_x}\int_{\mathbb{R}^2} |xv^N|^2|Y_N|\d x \leqslant \frac{1}{8}\left|v^N\right|_{\Sigma}^2+c_1|\e^{-2Y_N}|_{\mathbb{L}^\infty_x}^{2+p_1}|\e^{2Y_N}|_{\mathbb{L}^\infty_x}^{1+p_1} \left|Y_N\right|_{\W^{1-s,q}_x}^{p_1}|v_0|_{\mathbb{L}^2_x}^2.
        \end{equation*}
        Using \cref{Lemma_bound_from_moment,Lem:convergence-Y_N-W(1-s;q)} for the $\W^{1-s,q}_x$ norm and \cref{Lem:Conv-W(1-kappa;inf)-exp(aYN),Cor:moment-boud-expYN-norm-Linf} for the terms involving exponentials of $Y_N$, we obtain that $C_{1,N}=c_1|\e^{-2Y_N}|_{\mathbb{L}^\infty_x}^{2+p_1}|\e^{2Y_N}|_{\mathbb{L}^\infty_x}^{1+p_1} \left|Y_N\right|_{\W^{1-s,q}_x}^{p_1}$ verifies \hyperlink{property_star}{Property (*)} and
        \begin{equation}
            |\e^{-2Y_N}|_{\mathbb{L}^\infty_x}|\e^{2Y_N}|_{\mathbb{L}^\infty_x}\int_{\mathbb{R}^2} |xv^N(t)|^2|Y_N|\d x\leqslant \frac{1}{8}\left|v^N(t)\right|_{\Sigma}^2+C_{1,N}|v_0|_{\mathbb{L}^2_x}^2.\label{Equation_energy_bound_x2_YN}
        \end{equation}
    
    \item Let $\kappa\in(0,1)$, using duality and \cref{product_rule} we have
        \begin{align*}
            \left|\int_{\mathbb{R}^2} |v^N(t)|^2:|\nabla Y_N|^2:\e^{2Y_N}\d x\right|&\leqslant \left|:|\nabla Y_N|^2:\right|_{\W^{-\kappa,\infty}_x}\left||v^N(t)|^2 \e^{2Y_N}\right|_{\W^{\kappa,1}_x}\\ 
            &\lesssim  \left|:|\nabla Y_N|^2:\right|_{\W^{-\kappa,\infty}_x}|v^N(t)|_{\W^{\kappa,2}_x}^2 \left|\e^{2Y_N}\right|_{1+\W^{\kappa,\infty}_x}\\ 
            &\lesssim  \left|:|\nabla Y_N|^2:\right|_{\W^{-\kappa,\infty}_x}\left|\e^{2Y_N}\right|_{1+\W^{\kappa,\infty}_x}|v^N(t)|_{\mathbb{L}^{2}_x}^{2-2\kappa} |v^N(t)|_{\Sigma}^{2\kappa}.
        \end{align*}

    Let $\Tilde{C}_{2,N}=|\e^{-2Y_N}|_{\mathbb{L}^\infty} \left|:|\nabla Y_N|^2:\right|_{\W^{-\kappa,\infty}}\left|\e^{2Y_N}\right|_{1+\W^{\kappa,\infty}}$, by Young's inequality an mass conservation, there exists a constant $c_2>0$ such that
    \begin{equation*}
        |\e^{-2Y_N}|_{\mathbb{L}^\infty_x}\left|\int_{\mathbb{R}^2} |v^N(t)|^2:|\nabla Y_N|^2:\e^{2Y_N}\d x\right| \leqslant \frac{1}{8}|v^N(t)|_{\Sigma}^{2}+c_2 \Tilde{C}_{2,N}^{\frac{1}{1-\kappa}}|\e^{-2Y_N}|_{\mathbb{L}^\infty_x}|\e^{2Y_N}|_{\mathbb{L}^\infty_x}|v_0|_{\mathbb{L}^{2}_x}^{2}.
    \end{equation*}
    Let $C_{2,N}=c_2 \Tilde{C}_{2,N}^{\frac{1}{1-\kappa}} |\e^{2Y_N}|_{\mathbb{L}^\infty} |\e^{2Y_N}|_{\mathbb{L}^\infty}$, using \cref{Lem:Conv-W(1-kappa;inf)-exp(aYN),Cor:moment-boud-expYN-norm-Linf} for the terms involving exponentials of $Y_N$ and Sobolev embeddings, \cref{Cor:cauchy-criter-W(-kappa;infty),Lemma_bound_from_moment} for the renormalized term, it follows $C_{2,N}$ verifies \hyperlink{property_star}{Property (*)} and 
    \begin{equation}
         |\e^{-2Y_N}|_{\mathbb{L}^\infty_x}\left|\int_\mathbb{R}^2 |v^N(t)|^2:|\nabla Y_N|^2:\e^{2Y_N}\d x\right| \leqslant \frac{1}{8}|v^N(t)|_{\Sigma}^{2}+C_{2,N}|v_0|_{\mathbb{L}^{2}_x}^{2}. \label{Equation_energy_bound_renormalized}
    \end{equation}
\end{itemize}

Putting \cref{Equation_energy_bound_x2_YN,Equation_energy_bound_renormalized} back in \cref{Equation_bound_Sigma_from_E}, we obtain
$$|v^N(t)|_\Sigma^2 \leqslant 8|\e^{-2Y_N}|_{\mathbb{L}^\infty_x}\Tilde{\mathcal{E}}_N(v_0) + 8\left(C_{1,N}+C_{2,N}\right)|v_0|_{\mathbb{L}^2_x}^2.$$
As the right hand side is independent of $t$, for all $N\in\N$, $v^N$ is global and we have
\begin{equation}
    |v^N|_{\mathbb{L}^\infty_t\Sigma}^2 \leqslant 8|\e^{-2Y_N}|_{\mathbb{L}^\infty}\Tilde{\mathcal{E}}_N(v_0) + 8\left(C_{1,N}+C_{2,N}\right)|v_0|_{\mathbb{L}^2}^2. \label{Equation_energy_bound_defocusing_E(v0)}
\end{equation}

Now, recall the definition of the transform energy, 
\begin{equation}
    \Tilde{\mathcal{E}}_N(v_0) =\frac{1}{2}\displaystyle\int_{\mathbb{R}^2}\left( \left|\nabla v_0\right|^2+\left|xv_0\right|^2(1-Y_N)-:|\nabla Y_N|^2: \left|v_0\right|^2 \right)\e^{2Y_N} \d x\label{Eq:Energy-linear}
\end{equation}
so that using \cref{product_rule,product-rule-W(-k;r)-Wkp}, Sobolev embeddings and interpolation, for some arbitrary $\kappa\in(0,1)$, we obtain
\begin{equation}
    2\Tilde{\mathcal{E}}_N(v_0)\leqslant \left(1+|Y_N|_{\mathbb{L}^\infty_x}+\left|:|\nabla Y_N|^2:\right|_{\W^{-\kappa,\infty}_x}\right)|\e^{2Y_N}|_{1+\W^{\kappa,\infty}_x}|v_0|_{\Sigma}^{2}. \label{Equation_energy_bound_E(v0)_lin}
\end{equation}
Define $C_{N}=8\left(C_{1,N}+C_{2,N}\right)+ 4|\e^{-2Y_N}|_{\mathbb{L}^\infty_x} \left(1+|Y_N|_{\mathbb{L}^\infty_x}+\left|:|\nabla Y_N|^2:\right|_{\W^{-\kappa,\infty}_x}\right) |\e^{2Y_N}|_{1+\W^{\kappa,\infty}_x}$, then as previously, $C_{N}$ verifies \hyperlink{property_star}{Property (*)} and using \cref{Equation_energy_bound_defocusing_E(v0),Equation_energy_bound_E(v0)_lin}, we obtain$$|v^N|_{\mathbb{L}^\infty_t\Sigma}^2 \leqslant C_N|v_0|_{\Sigma}^2.$$
Taking the square root and $\mathbb{L}^{p}_\omega$ norm of both side and using Hölder's inequality, we finally obtain
$$\sup_{N\in\N}\left|v^N\right|_{\mathbb{L}^p_\omega\mathbb{L}^\infty_t\Sigma}\leqslant C |v_0|_{\mathbb{L}^{p_0}_\omega\Sigma}$$
using that all moments of $C_N$ are bounded uniformly in $N$ by \hyperlink{property_star}{Property (*)}.\endproof

\end{demo}

\begin{demo}{of \cref{Proposition_energy_bound_defocusing}.}

    As in the proof of \cref{Proposition_energy_bound_linear}, as $\lambda\leqslant 0$, for every $N\in\N$, $v^N$ is global and \cref{Equation_energy_bound_defocusing_E(v0)} holds. It remains to bound the term $\Tilde{\mathcal{E}}_N(v_0)$. Now recall that
    $$ \Tilde{\mathcal{E}}_N(v_0) = \Tilde{\mathcal{E}}_N^{lin}(v_0) - \frac{\lambda}{4}\int_{\mathbb{R}^2}\left|v\right|^{4} \e^{4 Y_N} \d x$$
    where $\Tilde{\mathcal{E}}_N^{lin}(v_0)$ is defined in \cref{Eq:Energy-linear} and verifies \cref{Equation_energy_bound_E(v0)_lin}. Using the same notation as in the proof of \cref{Proposition_energy_bound_linear}, we obtain
   \begin{equation}
       |v^N|_{\mathbb{L}^\infty_\R\Sigma}^2 \leqslant C_N|v_0|_{\Sigma}^2-2\lambda|\e^{-2Y_N}|_{\mathbb{L}^\infty_x}\int_{\mathbb{R}^2}\left|v_0\right|^{4} \e^{4 Y_N}\d x.
   \end{equation}
    Using \cref{Lemma_Gagliardo_Nirenberg}, we have
    \begin{equation}
        |v^N|_{\mathbb{L}^\infty_t\Sigma}^2 \leqslant C_N|v_0|_{\Sigma}^2-\lambda|\e^{-2Y_N}|_{\mathbb{L}^\infty_x}|\e^{4Y_N}|_{\mathbb{L}^\infty_x}|v_0|_{\Sigma}^4.\label{Equation_energy_bound_defocusing_final}
    \end{equation}
    As $C_N$ verifies \hyperlink{property_star}{Property (*)}, the claim follows from \cref{Equation_energy_bound_defocusing_final,Lem:Conv-W(1-kappa;inf)-exp(aYN)}.\endproof
\end{demo}

\begin{demo}{of \cref{Proposition_energy_bound_focusing}.}

    In the case where $\lambda>0$, \cref{Equation_bound_Sigma_from_E} does not hold, but instead, we have
    \begin{equation}
    \begin{split}
        \frac{1}{2}|v^N(t)|^2_\Sigma&\leqslant |\e^{-2Y_N}|_{\mathbb{L}^\infty}\Tilde{\mathcal{E}}_N(v_0)+|\e^{-2Y_N}|_{\mathbb{L}^\infty_x}|\e^{2Y_N}|_{\mathbb{L}^\infty_X}\int_{\mathbb{R}^2} |xv^N(t)|^2|Y_N|\d x\\ 
    &+\left||\e^{-2Y_N}|_{\mathbb{L}^\infty_x}\int_{\mathbb{R}^2} |v^N(t)|^2:|\nabla Y_N|^2:\e^{2Y_N}\d x\right|+\frac{\lambda|\e^{-2Y_N}|_{\mathbb{L}^\infty_x}}{4}\int_{\mathbb{R}^2} \left|v^N(t)\right|^4 \e^{4 Y_N}\d x.
    \end{split}\label{Equation_bound_Sigma_from_E_focusing}
    \end{equation}
    Using \cref{Lemma_Gagliardo_Nirenberg} and the mass conservation, we obtain    
    \begin{equation}
        \frac{\lambda|\e^{-2Y_N}|_{\mathbb{L}^\infty_x}}{4}\int_{\mathbb{R}^2} \left|v^N(t)\right|^4 \e^{4 Y_N} \d x \leqslant \frac{\lambda|\e^{-2Y_N}|_{\mathbb{L}^\infty_x}^2|\e^{4Y_N}|_{\mathbb{L}^\infty_x}}{8} \Tilde{M}_N(v_0)|v^N(t)|^2_\Sigma.\label{Equation_energy_bound_focusing_nonlinear}
    \end{equation}
    For $L>0$, let $\overline{\Omega}_L = \left\{\lambda L^2|\e^{-2Y}|_{\mathbb{L}^\infty_x}^2|\e^{2Y}|_{\mathbb{L}^\infty_x}|\e^{4Y}|_{\mathbb{L}^\infty_x}<4\right\}$. Then for any $\omega\in\overline{\Omega}_L$, there exists $N_0(\omega)\in\N$ and $\eps\in\left(0,\frac{1}{2}\right)$ such that for any $v_0\in\Sigma$ verifying $|v_0|_{\mathbb{L}^2_{x}}\leqslant L$, it holds
    $$\forall N\geqslant N_0(\omega),\; \frac{\lambda|\e^{-2Y_N}|_{\mathbb{L}^\infty_x}}{4}\int_{\mathbb{R}^2} \left|v^N(t)\right|^4 \e^{4 Y_N} \d x \leqslant \frac{\lambda|\e^{-2Y_N}|_{\mathbb{L}^\infty_x}^2|\e^{2Y_N}|_{\mathbb{L}^\infty_x}|\e^{4Y_N}|_{\mathbb{L}^\infty_x}}{8} |v_0|_{\mathbb{L}^2_{x}}^2\leqslant \frac{1}{2}-\eps.$$
    Then, \cref{Equation_bound_Sigma_from_E_focusing} gives on $\overline{\Omega}_L$, 
    \begin{equation}
    \begin{split}
        \forall N\geqslant N_0,\; \eps|v^N(t)|^2_\Sigma&\leqslant |\e^{-2Y_N}|_{\mathbb{L}^\infty_x}\Tilde{\mathcal{E}}_N(v_0)+|\e^{-2Y_N}|_{\mathbb{L}^\infty_x}|\e^{2Y_N}|_{\mathbb{L}^\infty_x}\int_{\mathbb{R}^2} |xv^N(t)|^2|Y_N|\d x\\ 
        &+\left||\e^{-2Y_N}|_{\mathbb{L}^\infty_x}\int_{\mathbb{R}^2} |v^N(t)|^2:|\nabla Y_N|^2:\e^{2Y_N}\d x\right|.
    \end{split}\label{Equation_bound_Sigma_from_E_focusing_linearRHS}
    \end{equation}
    A slight modification of the proof of \cref{Proposition_energy_bound_linear}, 
    shows there exists $\Tilde{C}_{1,N}$ and $\Tilde{C}_{2,N}$ verifying  \hyperlink{property_star}{Property (*)} and
    \begin{align}
        |\e^{-2Y_N}|_{\mathbb{L}^\infty_x}|\e^{2Y_N}|_{\mathbb{L}^\infty_x}\int_{\mathbb{R}^2} |xv^N(t)|^2|Y_N|\d x\leqslant \frac{\eps}{4}|v^N(t)|^2_\Sigma + \Tilde{C}_{1,N}|v_0|^2_{\mathbb{L}^2_{x}},\label{Equation_energy_bound_x2YN_focusing}\\ 
        |\e^{-2Y_N}|_{\mathbb{L}^\infty_x}\left|\int_{\mathbb{R}^2} |v^N(t)|^2:|\nabla Y_N|^2:\e^{2Y_N}\d x\right|\leqslant \frac{\eps}{4}|v^N(t)|^2_\Sigma + \Tilde{C}_{2,N}|v_0|^2_{\mathbb{L}^2_{x}}.\label{Equation_energy_bound_renormalized_focusing}
    \end{align}
    Moreover, in that case
    $$ \Tilde{\mathcal{E}}_N(v_0) = \Tilde{\mathcal{E}}_N^{lin}(v_0) - \frac{\lambda}{4}\int_{\mathbb{R}^2}\left|v\right|^{4} \e^{4 Y_N} \d x \leqslant \Tilde{\mathcal{E}}_N^{lin}(v_0)$$
    where $\Tilde{\mathcal{E}}_N^{lin}(v_0)$ is defined in \cref{Eq:Energy-linear} and verifies \cref{Equation_energy_bound_E(v0)_lin}. Define the random variable $\Tilde{C}_N = \frac{2}{\eps} \left(\left(1+|Y_N|_{\mathbb{L}^\infty_x}+\left|:|\nabla Y_N|^2:\right|_{\W^{-\kappa,\infty}_x}\right)|\e^{2Y_N}|_{1+\W^{\kappa,\infty}_x}+\Tilde{C}_{1,N}+\Tilde{C}_{2,N}\right)$ for some arbitrary $\kappa\in(0,1)$, it verifies \hyperlink{property_star}{Property (*)} as previously. Then, putting \cref{Equation_energy_bound_E(v0)_lin,Equation_energy_bound_x2YN_focusing,,Equation_energy_bound_renormalized_focusing} in \cref{Equation_bound_Sigma_from_E_focusing}, we obtain for every $\omega\in\overline{\Omega}_L$ and every $N\geqslant N_0(\omega)$,
   $$\forall t\in\mathbb{R},\; |v^N(\omega,t)|^2_\Sigma\leqslant \Tilde{C}_N(\omega) |v_0|_\Sigma^2.$$
   The claim follows as usual.
   \endproof

\end{demo}

\section{Proofs of Proposition \ref{Prop:diverging-bound-W^(sigma;2)-defocusing} and \ref{Proposition-Bound-L2-Diff-cubic-defocusing}}\label{appendix_proof_bounds_cubic}

\begin{lemme}\label{Lemma-diverging-bound-W^(sigma;2)}
    Let $a<0<b$, $\lambda\leqslant 0$, $\sigma\in\left(\frac{3}{2},2\right)$. For any $v_0\in\W^{2,2}$, there exists a positive random variable $C(v_0)$ almost surely finite such that for any $\kappa\in(0,1]$ there exists a positive random variable $C'(\kappa,v_0)$ almost surely finite, such that 
    $$\forall N\in\N,\; |v^N|_{\mathbb{L}^\infty_{t}\W^{\sigma,2}_x} \leqslant C'(\kappa,v_0) \lambda_N^{C(v_0)\kappa}.$$
    Moreover, $C$ and $C'(\kappa,\cdot)$ are nondecreasing functions of $|v_0|_{\W^{2,2}_x}$.
    
\end{lemme}

\begin{demo}{ of \cref{Lemma-diverging-bound-W^(sigma;2)}.}

    In what follows, we always assume conclusions of \cref{Proposition_energy_bound_defocusing} to hold. We will show a bound for $|v^N|_{\W^{\sigma,2}_x}$ uniform on $t\in[0,b]$, then reversing time, we obtain a bound of the same type uniformly in $t\in[a,0]$. Let $v_0\in\W^{2,2}$. Set $\mathfrak{C} = C(1+|v_0|_\Sigma)^2\geqslant \sup_{N\in\N} \left|v^N\right|_{\mathbb{L}^\infty_t\Sigma}$ the upper bound given by  \cref{Proposition_energy_bound_defocusing} and let $N\in\N$.\\ 

    Define $w^N=\p_t v^N$, then \cref{Eq:2D-transformed-renorm-regularized} gives
    \begin{align*}
        |v^N|_{\W^{\sigma,2}_x} &\leqslant |w^N|_{\mathbb{L}^2_{x}} + 2|\nabla v^N\cdot \nabla Y_N|_{\W^{\sigma-2,2}_x} + |x v^N\cdot x Y_N|_{\W^{\sigma-2,2}_x} \\ 
        &+ |v^N:|\nabla Y_N|^2:|_{\W^{\sigma-2,2}_x} + |\lambda| \left||v^N|^{2} v^N \e^{2 Y_N}\right|_{\mathbb{L}^2_{x}}.
    \end{align*}
    \Cref{Theorem-local-wellposedness,Proposition_energy_bound_defocusing} imply $w^N\in\mathcal{C}(\R,\mathbb{L}^2(\mathbb{R}^2))$, so we can define $W_N = \Tilde{M}_N(w^N)$. Then, by Young's inequality, we have  $$|w^N|_{\mathbb{L}^2_{x}}\leqslant \frac{W_N}{2} + \frac{\left|e^{-2Y_N}\right|_{\mathbb{L}^\infty_{x}}}{2}\leqslant \frac{W_N}{2} + \frac{\sup_{N\in\N}\left|e^{-2Y_N}\right|_{\mathbb{L}^\infty_{x}}}{2}.$$
    For the nonlinear term, by Sobolev embeddings, there exists a constant $c_0>0$ such that
    \begin{equation}
        \left||v^N|^{2} v^N \e^{2 Y_N}\right|_{\mathbb{L}^2_{x}} \leqslant c_0\left|e^{2 Y_N}\right|_{\mathbb{L}^\infty_{x}} \left| v^N\right|_\Sigma^{3}\leqslant c_0 \mathfrak{C}^{3} \sup_{N\in\N}\left|e^{2 Y_N}\right|_{\mathbb{L}^\infty_{x}}.\label{Proof-proposition-diverging-bound-estimate-nonlinear-term}
    \end{equation}
    Using \cref{product-rule-W(-k;r)-Wkp,Corollary_action_D_and_x_on_Wsp}, there exists a constant $c_1>0$ depending only of $\sigma$ such that:
    $$|\nabla v^N\cdot \nabla Y_N|_{\W^{\sigma-2,2}_x}\leqslant c_1 \left|v^N\right|_{\W^{3-\sigma,2}_x}\left|\nabla Y_N\right|_{\W^{\sigma-2,\infty}_x},$$
    $$|x v^N\cdot x Y_N|_{\W^{\sigma-2,2}_x}\leqslant c_1 \left|v^N\right|_{\W^{3-\sigma,2}_x}\left|x Y_N\right|_{\W^{\sigma-2,\infty}_x}$$
    and
    $$|v^N:|\nabla Y_N|^2:|_{\W^{\sigma-2,2}_x}\leqslant c_1 \left|v^N\right|_{\W^{2-\sigma,2}_x}\left|:|\nabla Y_N|^2:\right|_{\W^{\sigma-2,\infty}_x}.$$
    Let $\theta=\frac{2-\sigma}{1-\sigma}$. As $\sigma>\frac{3}{2}$ we have $0<2-\sigma<3-\sigma<\sigma$, $\theta\in(0,1)$ and $3-\sigma=\sigma\theta+(1-\theta)$. By interpolation, there exists $c_2>0$ such that
    \begin{align}
        |\nabla v^N\cdot \nabla Y_N|_{\W^{\sigma-2,2}_x}&\leqslant c_2\sup_{N\in\N}\left|\nabla Y_N\right|_{\W^{\sigma-2,\infty}_x}\mathfrak{C}^{1-\theta}|v^N|_{\W^{\sigma,2}_x}^\theta,\label{Proof-proposition-diverging-bound-estimate-grad-term}\\ 
        |x v^N\cdot x Y_N|_{\W^{\sigma-2,2}_x}&\leqslant c_2\sup_{N\in\N}\left|x Y_N\right|_{\W^{\sigma-2,\infty}_x}\mathfrak{C}^{1-\theta}|v^N|_{\W^{\sigma,2}_x}^\theta,\label{Proof-proposition-diverging-bound-estimate-x-term}\\ 
        |v^N:|\nabla Y_N|^2:|_{\W^{\sigma-2,2}_x}&\leqslant c_2 \sup_{N\in\N}\left|:|\nabla Y_N|^2:\right|_{\W^{\sigma-2,\infty}_x}\mathfrak{C}^{1-\theta}|v^N|_{\W^{\sigma,2}_x}^\theta.\label{Proof-proposition-diverging-bound-estimate-renormalized-term}
    \end{align}
    By \cref{Cor:convergence-Y_N-W(1-kappa;infty),Lem:Conv-W(1-kappa;inf)-exp(aYN),Cor:cauchy-criter-W(-kappa;infty)}, almost surely
    \begin{equation}
        \sup_{N\in\N}\max\left(\left|e^{-2Y_N}\right|_{\mathbb{L}^\infty_{x}},\left|e^{2 Y_N}\right|_{\mathbb{L}^\infty_{x}},\left|\nabla Y_N\right|_{\W^{\sigma-2,\infty}_x},\left|x Y_N\right|_{\W^{\sigma-2,\infty}_x},\left|:|\nabla Y_N|^2:\right|_{\W^{\sigma-2,\infty}_x}\right)<+\infty.\label{eq:condition-1-lemma-diverging-defocusing}
    \end{equation}
    By Young's inequality in \cref{Proof-proposition-diverging-bound-estimate-grad-term,Proof-proposition-diverging-bound-estimate-x-term,Proof-proposition-diverging-bound-estimate-renormalized-term}, and by \cref{Proof-proposition-diverging-bound-estimate-nonlinear-term}, there exists a positive random variable $C_0$, almost surely finite, such that
    $$|v^N|_{\W^{\sigma,2}_x} \leqslant \frac{W_N}{2} + \frac{|v^N|_{\W^{\sigma,2}_x}}{2} + \frac{C_0}{2}.$$
    Hence, we have
    $$|v^N|_{\W^{\sigma,2}_x} \leqslant W_N + C_0.$$
    Using that $w^N$ verify the equation:
    \begin{equation*}{~}
    \begin{cases}\mathrm{i}\p_t w^N&+ Hw^N +2\nabla Y_N\cdot \nabla w^N+ xY_N\cdot x w^N +:|\nabla Y_N|^2: w^N \\ 
    &+\lambda \left(|v^N|^{2}w^N+2\re\left(v^N\overline{w^N}\right)v^N\right)\e^{2 Y_N}=0\\ w^N(0)&=\mathrm{i}\left(Hv_0+2\nabla v_0\cdot \nabla Y_N+x v_0\cdot x Y_N +:|\nabla Y_N|^2: v_0 + \lambda |v_0|^{2} v_0 e^{2 Y_N}\right)\end{cases} 
    \end{equation*}
    we obtain
    $$\frac{1}{2} \frac{\d}{\d t} W_N\leqslant 2|\lambda| \int_{\mathbb{R}^2} \re\left(v^N\overline{w^N}\right)\im\left(v^N\overline{w^N}\right) e^{4 Y_N}\d x\leqslant 2|\lambda|  |v^N|^{2}_{\mathbb{L}^\infty_{x}}|\e^{2 Y_N}|_{\mathbb{L}^\infty_{x}} W_N.$$
    Using \cref{Generalized_Brezis-Gallouet_inequality}, there exists a positive random variable $C_1$, almost surely finite, such that
    $$\frac{\d}{\d t} W_N \leqslant C_1 \left(1+\ln\left(1+|v^N|_{\W^{\sigma,2}_x}\right)\right) W_N\leqslant C_1 \left(1+\ln\left(1 + W_N + C_0\right)\right) W_N.$$
    Define $G_N = 1 + \ln\left(1 + W_N + C_0\right)$, then it verifies almost surely $G_N' \leqslant C_1 G_N.$ Hence, Gronwall's lemma gives
    $$\forall t\in[0,b],\; G_N(t) \leqslant G_N(0)\e^{t C_1} \leqslant (1+\ln(1+C_0+W_N(0))\e^{b C_1}$$
    which implies
    $$|v^N|_{\mathbb{L}^\infty_{t}\W^{\sigma,2}_x}\leqslant \e^{\e^{b C_1}} (1+C_0+W_N(0))^{e^{b C_1}}.$$
    Now, let $\kappa\in(0,1)$. Remark that
    \begin{alignat*}{2}
        &W_N(0) &\leqslant |\e^{2 Y_N}|_{\mathbb{L}^\infty_{x}} |w_0^N|_{\mathbb{L}^2_{x}}^2 &\leqslant \sup_{n\in\N}|\e^{2 Y_n}|_{\mathbb{L}^\infty_{x}}\left(|v_0|_{\W^{2,2}_x}+|\lambda||v_0|_{\mathbb{L}^6_x}^3\sup_{n\in\N}|\e^{2 Y_n}|_{\mathbb{L}^\infty_{x}}\right.\\ 
        &&&+\left.\sup_{n\in\N}\left[\lambda_n^{-\kappa}\max\left(|\nabla v_0\cdot \nabla Y_n|_{\mathbb{L}^2_{x}},|x v_0\cdot x Y_n|_{\mathbb{L}^2_{x}}, |v_0:|\nabla Y_n|^2:|_{\mathbb{L}^2_{x}}\right)\right]\lambda_N^{\kappa}\right)^2.
    \end{alignat*}
    By \cref{Cor:estimate-norm-Lq-X-Y_N,estim-moment-norm-Lq-X.produit-wick}, almost surely
    $$\sup_{n\in\N}\left[\lambda_n^{-\kappa}\max\left(|\nabla v_0\cdot \nabla Y_n|_{\mathbb{L}^2_{x}},|x v_0\cdot x Y_n|_{\mathbb{L}^2_{x}}, |v_0:|\nabla Y_n|^2:|_{\mathbb{L}^2_{x}}\right)\right]<+\infty.$$
    Moreover, there exists a positive random variable $C'_0(\kappa)$, depending on $v_0$ through $|v_0|_{\W^{2,2}_x}$ non-decreasingly, almost surely finite, such that $W_N(0) \leqslant C'_0(\kappa)\left(1+\lambda_N^{2\kappa}\right),$ and as $\lambda_N\geqslant 1$, we obtain 
    $$|v^N|_{\mathbb{L}^\infty_{t}\W^{\sigma,2}_x}\leqslant (2\e(1+C_0+C'_0(\kappa)))^{\e^{b C_1}} \lambda_N^{2\kappa e^{b C_1}}.$$
    Setting $C(v_0)=2e^{b C_1}$ and $C'(\kappa,v_0)=(2\e(1+C_0+C'_0(\kappa)))^{\e^{b C_1}}$, we obtain the claim. \endproof
\end{demo}

\begin{demo}{of \cref{Prop:diverging-bound-W^(sigma;2)-defocusing}.}

    For $n\in\N$, define $\sigma_n = 2-2^{-(n+2)}$. Hence, for any $a<0<b$ and any $\sigma\in(1,2)$ there exists $n\in\N$ such that $\mathbb{L}^\infty([-n,n],\W^{\sigma_n,2})$ is continuously embedded in $\mathbb{L}^\infty([a,b],\W^{\sigma,2})$. Thus it is sufficient to prove the claim for fixed $a<0<b$ and $\sigma\in\left(\frac{3}{2},2\right)$. Moreover, as $\W^{2,2} = \bigcup_{n\in\N} B_{\W^{2,2}}(0,n)$, it is sufficient to prove the claim uniformly in $v_0\in B_{\W^{2,2}}(0,L)$ a $\W^{2,2}$-ball of fixed radius $L\in\N$.\\ 

    Let $L\in\N$. For $j\in\N$, let $\kappa_j = 2^{-j}$. As $\lambda_N>1$ for all $N$, it is sufficient to show conclusion only for $(\kappa_j)_{j\in\N^*}$. Recall the random variable $C(v_0)$ and $C'(\kappa,v_0)$ defined in \cref{Lemma-diverging-bound-W^(sigma;2)}. Denote by $C_L$ and $C'_L(\kappa)$ the value of $C(v'_0)$ and $C'(\kappa,v'_0)$ respectively for a fixed $v'_0\in \W^{2,2}$ of $\W^{2,2}_x$-norm equal to $L$, recall that $C$ and $C'$ are nondecreasing functions of $|v_0|_{\W^{2,2}_x}$. For $k\in\N$, let $\Tilde{\Omega}_k = \left\{2^k\leqslant C_L<2^{k+1}\right\}$ and let $\Tilde{\Omega}_{-1}=\{C_L<1\}$. By definition of $C_L$, $\Tilde{\Omega}=\bigcup_{k\geqslant-1}\Tilde{\Omega}_k$ is of probability one. Let $\kappa'\in(0,1)$ and $k\geqslant -1$, in view of \cref{Lemma-diverging-bound-W^(sigma;2)}, on $\Tilde{\Omega}_k$, it holds
    $$\forall N\in\N,\; \forall v_0\in B_{\W^{2,2}}(0,L),\; \left|v^N\right|_{\mathbb{L}^\infty_{t}\W^{\sigma,2}_x} \leqslant C'_L(\kappa')\lambda_N^{C_L\kappa'}\leqslant C'_L(\kappa')\lambda_N^{2^{k+1}\kappa'}.$$
    Let $\Omega'_\kappa = \{C'_L(\kappa)<+\infty\}$, from \cref{Lemma-diverging-bound-W^(sigma;2)} it is of probability one for any $\kappa\in(0,1]$. Thus $\Omega' = \bigcap_{n\in\N} \Omega'_{2^{-n}}$ is of probability one too. On $\Omega'\cap\Tilde{\Omega}_k$, it holds
    $$\forall N\in\N,\; \forall v_0\in B_{\W^{2,2}}(0,L),\; \left|v^N\right|_{\mathbb{L}^\infty_{t}\W^{\sigma,2}_x} \leqslant C'_L(2^{-(k+j+1)})\lambda_N^{\kappa_j}.$$
    Hence, on $\Tilde{\Omega}\cap\Omega'$, which is of probability one, it holds
    $$\forall j\in\N,\;\forall N\in\N,\; \forall v_0\in B_{\W^{2,2}}(0,L),\; \left|v^N\right|_{\mathbb{L}^\infty_{t}\W^{\sigma,2}_x} \leqslant C_L(\kappa_j)\lambda_N^{\kappa_j},$$
    for a certain positive random variable $C_L(\kappa_j)$ finite on $\Tilde{\Omega}\cap\Omega'$, which allows to conclude.\endproof
    
\end{demo}

\begin{rem}
    If one tries to apply the proof of \cref{Lemma-diverging-bound-W^(sigma;2)} in the case $\sigma = 2$, it would end with a random variable depending of $\kappa$ in the exponent, thus the proof of \cref{Prop:diverging-bound-W^(sigma;2)-defocusing} would not hold.
\end{rem}

\begin{lemme}\label{Lemma-Bound-L2-Diff-cubic}

    Let $a<0<b$ and $\lambda\leqslant0$. For any $v_0\in\W^{2,2}$, there exists a positive random variable $C(v_0)$ almost surely finite such that for any $\kappa\in(0,1]$ and any $\delta\in\left(0,\frac{1}{3}\right)$ there exists a positive random variable $C''(\kappa,\delta,v_0)$ almost surely finite, such that 
    $$\forall N>M\geqslant 0,\; |v^N-v^M|_{\mathbb{L}^\infty_{t}\mathbb{L}^2_{x}} \leqslant C''(\kappa,\delta,v_0) \lambda_N^{C(v_0)\kappa}\lambda_M^{-\delta}.$$
    Moreover, $C$ and $C''$ are nondecreasing functions of $|v_0|_{\W^{2,2}_x}$.
    
\end{lemme}

\begin{demo}{ of \cref{Lemma-Bound-L2-Diff-cubic}.}

    Assume in what follows that both conclusions of \cref{Proposition_energy_bound_defocusing,Prop:diverging-bound-W^(sigma;2)-defocusing} hold. Let $\delta\in\left(0,\frac{1}{3}\right)$, $\delta'\in\left(2\delta,\frac{2}{3}\right)$, $\kappa'\in\left(\frac{3\delta'}{2},1\right)$ and $\sigma\in(1+\kappa',2)$. Define $\theta = \frac{1+\kappa'}{\sigma}\in(0,1)$. Let $N>M\geqslant 2$. As in the proof of \cref{Lemma-diverging-bound-W^(sigma;2)}, let $v_0\in\W^{2,2}$ and set $\mathfrak{C} = C(1+|v_0|_\Sigma)^2\geqslant \sup_{n\in\N} \left|v^n\right|_{\mathbb{L}^\infty_t\Sigma}$ the upper bound given by  \cref{Proposition_energy_bound_defocusing}. Define $R= v^N-v^M$, then it verifies 
    \begin{equation}
    \begin{cases}-i\p_t R =& HR+2\nabla Y_N\cdot \nabla R+2\nabla (Y_N-Y_M)\cdot \nabla v^M+ xY_N\cdot x R +:|\nabla Y_N|^2: R \\ 
    &+ \lambda\left|v^N\right|^{2} R \e^{2 Y_N} + x(Y_N-Y_M)\cdot x v^M + \left(:|\nabla Y_N|^2:-:|\nabla Y_M|^2: \right)v^M\\  
    &+\lambda \left(\left|v^N\right|^{2}-\left|v^M\right|^{2}\right)v^M\e^{2 Y_N} + \lambda \left|v^M\right|^{2}v^M\left(\e^{2 Y_N}-\e^{2 Y_M}\right)
    \\ R(0)=0.\end{cases} \label{2d_transf_renorm_regu_eq_non-lin-diff-N-M}
    \end{equation}
    Now, set $W_N = \int_{\mathbb{R}^2} |R|^2 \e^{2Y_N}\d x$. Taking the time derivative, we obtain 
    \begin{equation}\label{Eq:time-derivative-proof-Cauchy-cubic-defocusing}
        \begin{split}
            \frac{1}{2} W_N' &\leqslant 2 |\nabla(Y_N-Y_M)|_{\W^{-\kappa',\infty}}\left|\overline{R}\nabla v^M\e^{2Y_N}\right|_{\W^{\kappa',1}} + |x(Y_N-Y_M)|_{\W^{-\kappa',\infty}}\left|x \overline{R} v^M\e^{2Y_N}\right|_{\W^{\kappa',1}} \\ 
            &+ \left|:|\nabla Y_N|^2:-:|\nabla Y_M|^2:\right|_{\W^{-\kappa',\infty}}\left|\overline{R} v^M\e^{2Y_N}\right|_{\W^{\kappa',1}} \\ 
            &+ \left|\lambda\int_{\mathbb{R}^2} \left(\left|v^N\right|^{2}-\left|v^M\right|^{2}\right)v^M\overline{R}\e^{4Y_N}\d x\right| + \left|\lambda\int_{\mathbb{R}^2} \left|v^M\right|^{2}v^M\left(\e^{2 Y_N}-\e^{2 Y_M}\right)\overline{R}\e^{2Y_N}\d x\right|.
        \end{split}
    \end{equation}
    First, by Cauchy-Schwarz inequality and Sobolev embeddings, there exists $c_0>0$ such that
    \begin{equation}
        \left|\lambda\int_{\mathbb{R}^2} \left|v^M\right|^{2}v^M\left(\e^{2 Y_N}-\e^{2 Y_M}\right)\overline{R}\e^{2Y_N}\d x\right| \leqslant c_0 \mathfrak{C}^4 \sup_{n\in\N} \left|e^{2 Y_n}\right|_{\mathbb{L}^\infty_{x}} \sup_{n\in\N} \lambda_n^{\delta'}\left|e^{2 Y_n}-e^{2 Y}\right|_{\mathbb{L}^\infty_{x}} \lambda_M^{-\delta'}.\label{eq:diff-exp-cauchy-estimate}
    \end{equation}
    Using \cref{product_rule,action_D&x/Wsp}, by interpolation, there exists $c_1>0$ such that
    \begin{align}
        \max\left(\left|\overline{R}\nabla v^M\e^{2Y_N}\right|_{\W^{\kappa',1}_x},\left|x\overline{R} v^M\e^{2Y_N}\right|_{\W^{\kappa',1}_x}\right)&\leqslant c_1 \sup_{n\in\N} \left|e^{2 Y_n}\right|_{1+\W^{\kappa',\infty}_x} \mathfrak{C}^{2-\theta} \left|v^M\right|_{\mathbb{L}^\infty_t\W^{\sigma,2}_x}^\theta\label{eq:Rgrad-Rx-cauchy-estimate}\\ 
        \left|\overline{R} v^M\e^{2Y_N}\right|_{\W^{\kappa',1}_x}&\leqslant c_1 \sup_{n\in\N} \left|e^{2 Y_n}\right|_{1+\W^{\kappa',\infty}_x} \mathfrak{C}^{2}.\label{eq:Rv-cauchy-estimate}
    \end{align}
    Moreover, we have 
    \begin{align}
        |\nabla(Y_N-Y_M)|_{\W^{-\kappa',\infty}_x} &\leqslant 2 \sup_{n\in\N} \lambda_n^{\delta'} |\nabla(Y_n-Y)|_{\W^{-\kappa',\infty}_x} \lambda_M^{-\delta'}\label{eq:grad-proof-cauchy-estimate}\\ 
        |x(Y_N-Y_M)|_{\W^{-\kappa',\infty}_x}&\leqslant 2 \sup_{n\in\N} \lambda_n^{\delta'} |x(Y_n-Y)|_{\W^{-\kappa',\infty}_x} \lambda_M^{-\delta'}\label{eq:x-proof-cauchy-estimate}\\ 
        \left|:|\nabla Y_N|^2:-:|\nabla Y_M|^2:\right|_{\W^{-\kappa',\infty}_x}&\leqslant 2 \sup_{n\in\N} \lambda_n^{\delta'} \left|:|\nabla Y_n|^2:-:|\nabla Y|^2:\right|_{\W^{-\kappa',\infty}_x} \lambda_M^{-\delta'}\label{eq:Wick-proof-cauchy-estimate}\\ 
        \left|\int_{\mathbb{R}^2} \left(\left|v^N\right|^{2}-\left|v^M\right|^{2}\right)v^M\overline{R}\e^{4Y_N}\d x\right| &\leqslant \sup_{n\in\N} \left|e^{2 Y_n}\right|_{\mathbb{L}^\infty_{x}}\left|v^M\right|_{\mathbb{L}^\infty_{x}}\left(\left|v^M\right|_{\mathbb{L}^\infty_{x}}+\left|v^N\right|_{\mathbb{L}^\infty_{x}}\right) W_N.\label{eq:Nonlin-proof-cauchy-estimate}
    \end{align}
    From \cref{Cor:convergence-Y_N-W(1-kappa;infty),Lem:Conv-W(1-kappa;inf)-exp(aYN),Cor:cauchy-criter-W(-kappa;infty)}, it follows almost surely
    $$\sup_{n\in\N} \lambda_n^{\delta'}\max\left(|\nabla(Y_n-Y)|_{\W^{-\kappa',\infty}_x},|x(Y_n-Y)|_{\W^{-\kappa',\infty}_x},\left|:|\nabla Y_n|^2:-:|\nabla Y|^2:\right|_{\W^{-\kappa',\infty}_x}\right)<+\infty$$
    and
    $$\sup_{n\in\N} \max\left(\lambda_n^{\delta'}\left|e^{2 Y_n}-e^{2 Y}\right|_{\mathbb{L}^\infty_{x}},\left|e^{2 Y_n}\right|_{1+\W^{\kappa',\infty}_x}\right)<+\infty.$$
    Thanks to \cref{Generalized_Brezis-Gallouet_inequality,Proposition_energy_bound_defocusing}, there exists a deterministic $c_2>0$ such that, almost surely, it holds
    $$\forall n\in\N,\; \left|v^n\right|_{\mathbb{L}^\infty_{x}}\leqslant c_2 \sqrt{1+\mathfrak{C}\left(1+\log\left(1+\left|v^n\right|_{\mathbb{L}^\infty_{t}\W^{\sigma,2}_x}\right)\right)}.$$
    Using \cref{Proposition_energy_bound_defocusing,Prop:diverging-bound-W^(sigma;2)-defocusing,eq:grad-proof-cauchy-estimate,Eq:time-derivative-proof-Cauchy-cubic-defocusing,eq:x-proof-cauchy-estimate,eq:Wick-proof-cauchy-estimate,eq:Nonlin-proof-cauchy-estimate,eq:diff-exp-cauchy-estimate,eq:Rgrad-Rx-cauchy-estimate,eq:Rv-cauchy-estimate}, there exists a positive random variable $C_0=C_0(v_0)$, almost surely finite, which is independent of $N$, $M$ and $t$, such that, almost surely, for any $\kappa\in(0,1)$,
    $$W_N' \leqslant C_0\left(1+\left|v^M\right|^{\theta}_{\mathbb{L}^\infty_{t}\W^{\sigma,2}_x}\right)\lambda_M^{-\delta'} + C_1 \left(1+\ln\left(1+C'_\kappa\lambda_N^\kappa\right)\right) W_N,$$
    where
    $$C_1 = 2 c_2^2 \sup_{n\in\N} \left|e^{2 Y_n}\right|_{\mathbb{L}^\infty_{x}} (1+\mathfrak{C})$$
    and $C'_\kappa$, given by \cref{Prop:diverging-bound-W^(sigma;2)-defocusing}, are almost surely finite random variables. Moreover, $C_0$, $C_1$ and $C'_\kappa$ are nondecreasing functions of $|v_0|_{\W^{2,2}_x}$. We now apply \cref{Prop:diverging-bound-W^(sigma;2)-defocusing} with an arbitrary $\kappa''\in\left(0,\min\left(1,\frac{\delta'-2\delta}{\theta}\right)\right)$ to obtain
    $$W_N' \leqslant \Tilde{C}_{\delta}\lambda_M^{-2\delta} + C_1 \left(1+\ln\left(1+C'_\kappa\lambda_N^\kappa\right)\right) W_N,$$
    where $\Tilde{C}_{\delta}$ is a positive random variable, almost surely finite, which is independent of $N$, $M$ and $t$, and is a nondecreasing function of $|v_0|_{\W^{2,2}_x}$. Then by Gronwall's lemma and \cref{Lem:Conv-W(1-kappa;inf)-exp(aYN)}, we conclude that there exists a positive random variable $C''(\kappa,\delta,v_0)$ almost surely finite, independent of $N$ and $M$ such that
    $$\left|v^N-v^M\right|_{\mathbb{L}^\infty_{t}\mathbb{L}^2_{x}}^2 \leqslant C''(\kappa,\delta,v_0)\lambda_N^{C_1\kappa}\lambda_M^{-\delta}.$$
    Moreover, $C''(\kappa,\delta,v_0)$ is a nondecreasing function of $|v_0|_{\W^{2,2}_x}$.\endproof

\end{demo}

\begin{demo}{ of \cref{Proposition-Bound-L2-Diff-cubic-defocusing}.}

    The proof is analogous to the proof of \cref{Prop:diverging-bound-W^(sigma;2)-defocusing}. \endproof
    
\end{demo}

\section*{Acknowledgments}

This work was supported by a public grant from the Fondation Mathématique Jacques Hadamard. Moreover, the author warmly thanks his PhD advisor, Anne de Bouard, for her constant support and helpful advices during the redaction of this paper.

\sloppy
\printbibliography

\end{document}